\documentclass[reqno,11pt]{amsart}
\usepackage{graphicx}
\usepackage{amscd}
\usepackage{amsmath}
\usepackage{caption}
\usepackage{amsfonts}
\usepackage{amssymb}
\usepackage{mathrsfs}
\usepackage{multicol}
\usepackage{color}

\newcommand{\be}{\begin{equation}}
\newcommand{\ee}{\end{equation}}
\newcommand{\ben}{\begin{eqnarray*}}
\newcommand{\een}{\end{eqnarray*}}

\newtheorem{theorem}{Theorem}
\newtheorem{lemma}{Lemma}
\newtheorem{remark}{Remark}
\newtheorem{corollary}{Corollary}
\newtheorem{proposition}{Proposition}
\newtheorem{definition}{Definition}
\allowdisplaybreaks
\definecolor{darkgreen}{rgb}{0.09, 0.45, 0.27}
\definecolor{debianred}{rgb}{0.84, 0.04, 0.33}
\allowdisplaybreaks[1]

\numberwithin{equation}{section}
\numberwithin{theorem}{section}
\begin{document}

%\begin{frontmatter}
%\hfill{\today}\bigskip 

\title[Positive radial solutions of a quasilinear problem with diffusion]
{Existence and nonexistence of positive radial solutions of a quasilinear Dirichlet problem\\ with diffusion}

%\thanks{This research was supported by the Italian MIUR project titled  {\it ``Metodi Va\-ria\-zio\-na\-li ed Equazioni Differenziali non Lineari''}.}

\author[L. Baldelli]{Laura Baldelli} \email{lbaldelli@impan.pl}

\author[V. Brizi]{Valentina Brizi} \email{valentina\_brizi@libero.it}

\author[R. Filippucci]{Roberta Filippucci} \email{roberta.filippucci@unipg.it}

\address[Filippucci]{Department of Mathematics -- University of Perugia --
Via Vanvitelli 1 -- 06123 Perugia,  Italy}
\address[Baldelli]{Institute of Mathematics -- Polish Academy of Sciences --
ul. Sniadeckich 8 -- 00-656 Warsaw, Poland}
\address[Brizi]{Department of Mathematics -- University of Perugia --
Via Vanvitelli 1 -- 06123 Perugia,  Italy}

\begin{abstract}
In this paper existence and nonexistence results of positive radial solutions of a Dirichlet $m$-Laplacian problem with different weights and 
a diffusion term inside the divergence  of the form $\big(a(|x|)+g(u)\big)^{-\gamma}$,  with $\gamma>0$ and $a$, $g$ positive functions satisfying natural growth conditions, are proved. Precisely, we obtain a new critical exponent $m^*_{\alpha,\beta,\gamma}$, which extends the one relative to case with no diffusion and it divides existence from nonexistence of positive radial solutions. The results are obtained via several tools such as a suitable modification of the celebrated blow up technique, Liouville type theorems, a fixed point theorem and a Poho\v zaev-Pucci-Serrin type identity. 
\end{abstract}


\maketitle

\section{Introduction}\label{intro}

In this paper we study existence and nonexistence of positive radial solutions of the 
following nonlinear elliptic problem
\begin{equation} \label{p1.1}
\begin{cases}
      -\mbox{div}\big(A(x,u)|\nabla u|^{m-2}\nabla u\big)=b(x,u) \quad &\text{in}\quad B_R,\\
      u=0 &\text{on} \quad\partial B_R,
\end{cases}\end{equation}
where $B_R$ denotes the open ball of radius $R>0$ centered at the origin in $\mathbb R^N$, $N>m$, $b$ is a continuous function and the differential operator involved is the $m$-Laplacian, namely $\Delta_m u=\mbox{div}(|\nabla u|^{m-2}\nabla u)$,
$m>1$, while $\nabla u$ denotes the gradient of $u$. 

Taking inspiration from \cite{Ba}, we are interested in diffusion terms $A$ of the form 
$$A(x,u)=\frac{|x|^\alpha}{\big(a(|x|)+g(u)\big)^\gamma}\quad\text{in}\quad\mathbb R^N\setminus\{0\}\times \mathbb R^+_0,$$
where  $\alpha \in\mathbb R$,  $\gamma>0$, with $g$ and $a$  continuous nonnegative functions satisfying suitable properties.
The main prototype for the diffusion term $A$ widely studied in literature is 
$A=(1+|u|)^{-\gamma}$, we refer to \cite{RHH,PSma} for a detailed discussion where measure data are taken under consideration.

Precisely, here we study positive radial solutions of the following Dirichlet problem
\begin{equation}\label{PROB}
\begin{cases}
  -\mbox{div}\biggl(\dfrac{|x|^\alpha|\nabla u|^{m-2}\nabla u}{\big(a(|x|)+g(u)\big)^\gamma} \biggr)=|x|^\beta u^p \quad &\text{in}\; B_R\setminus\{0\},
  \\ u=0 &\text{on}\; \partial B_R,
\end{cases}
\end{equation}
where $p>1$,  $\alpha,\beta\in\mathbb R$ and $\gamma\in(0,m-1)$.
We will say that $u$ is a solution of \eqref{PROB} if $u\in C^{1}(B_R)\cap C({\overline{B_R}})$, $u\ge0$ and $u$ solves the equation in problem \eqref{PROB} in the weak sense. Actually, we restrict our attention to  positive radial solutions $u$ of \eqref{PROB}, that is to positive functions $v=v(r)=u(|x|)$ such that  \begin{enumerate}
    \item[$(i)_s$] $v\in C^{1}[0,R)\cap C[0,R]$, 
    \item[$(ii)_s$] $r^{N+\alpha-1}|v'(r)|^{m-2}v'(r)\big[a(r)+g(v(r))\big]^{-\gamma}\in C^1(0,R)$
\end{enumerate}  and $v$ satisfies
\begin{equation} \label{probrad}
 \begin{cases}
  -\left(\dfrac{r^{N+\alpha-1}|v'(r)|^{m-2}v'(r)}{\big(a(r)+g(v(r))\big)^\gamma}\right)'=r^{N+\beta-1}v^p(r), \qquad 0< r<R
  \\v'(0)=0, \qquad v(R)=0.
   \end{cases}
\end{equation}
For the main properties of solutions of \eqref{probrad}, we refer to Section \ref{cap1}, where we also discuss the validity of condition $v'(0)=0$.

Differently from the case $\gamma=0$, as noted in \cite{Ba},  problem \eqref{PROB} cannot be attached with variational methods when $\gamma\neq0$ since the operator $\mbox{div}(A(x,u)|\nabla u|^{m-2}\nabla u)$ is well defined in $W^{1,m}(\Omega)$, but it may fail to be coercive on the same space when $u$ is large because of the properties assumed on $g$. Due to the lack of coercivity, the classical theory for elliptic operators acting between spaces in duality cannot be applied.
To overcome this difficulty, a blow up type technique turns out to be crucial in obtaining existence.  It is well known that the blow up method, due to Gidas and Spruck in the celebrated paper \cite{GS}, is based on producing a priori uniform bounds for positive solutions of Dirichlet problem in bounded domains, which are somehow equivalent to the validity of Liouville type results, namely nonexistence results in the entire space $\mathbb R^N$. In turn, existence follows by an application of degree theory based on the fixed point theorem by Krasnosel'skii.

Problem \eqref{PROB} in the semilinear case, i.e. $m=2$, without diffusion and weights, that is  $\gamma=0$ and $\alpha=\beta=0$ respectively, in bounded domains has been widely studied in huge of papers starting for the milestones papers  by Brezis, Nirenberg \cite{BN},  Trudinger in \cite{T68} and Aubin \cite{A76}, for  critical nonlinearities and by Ni and Serrin in \cite{NS} where they study also ground state positive solutions of $-\Delta u=f(u)$ for different types of nonlinearities $f$. 

The first Liouville  result for power type nonlinearities in the semilinear case was proved by Gidas and Spruck in \cite{GS} and states that 
{\it if $1< p<(N+2)/(N-2)=2^*-1$, $N>2$, where $2^*=2N/(N-2)$ is the Sobolev exponent for the Laplacian, then every nonnegative solution $u$ of class $C^2(\mathbb R^N)$ of the Lane-Emden equation $-\Delta u=u^p$ in $\mathbb R^N$ is such that $u\equiv0$}. The result is sharp. 

Later, the quasilinear case, $1<m<N$ and $\gamma=0$, which arises in many nonlinear phenomena such as in the theory of quasi-regular and quasi-conformal mappings, as well as a mathematical model of non-Newtonian fluids ($1<m<2$ pseudo plastic fluids while  $m>2$ dilating fluids such as blood), attracted much attention.  

Concerning the $m$-Laplacian problem $-\Delta_m u=u^p$ in $\mathbb R^N$, Mitidieri in \cite{mdokl98, mdokl} and Serrin and Zou in \cite{SZ} extended the above Liouville result obtaining nonexistence of nonnegative solutions if and only if  $$m-1< p<\frac{N+m}{N-m}=m^*-1,
\qquad N>m,\qquad  m^*=\frac{mN}{N-m},$$ 
where $m^*$ is the critical Sobolev exponent for the $m$-Laplacian, see Corollary II in \cite{SZ}.
On the other hand, Liouville type results for the  inequality $-\Delta_m u\ge u^p$ in $\mathbb R^N$ or in exterior domains involve no more the  Sobolev exponent $m^*$, but the Serrin exponent $m_*$
for the $m$-Laplacian case, that is 
$$m_*=\frac{m(N-1)}{N-m}(<m^*),$$
indeed, nonexistence for the inequality holds for $p\le m_*-1$.
For details, we refer to Mitidieri Poho\v zaev \cite{mdokl98},  \cite{mdokl},  to
Bidaut-V\'eron and Poho\v zaev \cite{B}, \cite{BP}  and to Serrin and Zou in \cite{SZ}.
Earlier nonexistence results for radially symmetric solutions have been established in \cite{NS}, while for Liouville type theorems for stable solutions or for solutions stable outside a compact set we refer to \cite{DFSV}.  In \cite{AMM}  a detailed analysis of the asymptotic behaviour of positive groundstate solutions with a nonlinearity subcritcal, critical and supercritical.

For $m$-Laplacian critical problems with different weights, we refer to the paper by  Clem\'ent, de Figueiredo and Mitidieri in \cite{CDM}, where they obtain the critical exponent associated to \eqref{probrad} with $\gamma=0$ given by
\begin{equation}\label{p*CDM}
m^*_{\alpha,\beta}:=\frac{m(N+\beta)}{N+\alpha-m}, \quad\qquad m-N<\alpha<\beta+1.
\end{equation}
We observe that $m^*_{\alpha,\beta}$ reduces to $m^*$ when $\alpha=\beta=0$, furthermore
$m^*_{\alpha,\beta}>m$ since $\beta-\alpha+m>0$. 

In particular, in  \cite{CCM98, CM}, the Authors prove non existence of positive radial solutions of 
\begin{equation}\label{liouv1}
  -\mbox{div}\bigl(|x|^\alpha|\nabla u|^{m-2}\nabla u\bigr)=|x|^\beta u^p \quad \text{in } \mathbb R^N
\end{equation}
under the assumption
\begin{enumerate}
\item [$(H_1)'$]$\quad m -1<p<m^*_{\alpha,\beta}-1,$
\end{enumerate}
while in Theorem 4.1 in \cite{CDM},   perform a detailed analysis of solutions of \eqref{liouv1} is performed yielding a nonexistence result in bounded domains for the critical case $p=m^*_{\alpha,\beta}-1$, by using variational techniques. We mention also  the paper by Egnell \cite{Eg} for additional results in this direction. 

Further existence and nonexistence results of positive solutions in the subcritical case for the Hardy–Hénon equation with Dirichlet boundary condition and with weights can be found in the recent paper \cite{CWZ}.

When a diffusion is introduced in the equation,  for the semilinear case of \eqref{PROB} in a generic bounded domain $\Omega$ and without weights,  that is $m=2$, $0<\gamma<1$ and $\alpha=\beta=0$, we mainly refer to \cite{BDO} where measure data are involved. To overcome the fact that the differential operator $\mbox{div}(\nabla u/(1+|u|)^\gamma)$ is not coercive on $H^1_0(\Omega)$ when $u$ is large,  see \cite{Por} for an explicit proof of this fact, the authors in \cite{BDO} need to work with approximate nondegenerate problems, using some a priori estimates since classical methods cannot be applied even if the datum $f$ is very regular. Further results, dealing with the same problem but with different order of summability for $f$, were proved  by Boccardo in \cite{Bocc} and by Boccardo and Brezis in \cite{BB} considering  the nonlinearity $f$ belonging to $L^\sigma(\Omega)$ with respectively $1\le\sigma\le N/2$ and $\sigma>N/2$. 

Then, in 2003, Alvino, Boccardo, Ferone, Orsina and Trombetti in \cite{BFR} studied the existence of positive solutions of the following problem 
\begin{equation}\label{pb_proto}
\begin{cases}
    -\mbox{div}\biggl(\dfrac{|\nabla u|^{m-2}\nabla u}{(1+|u|)^{\theta(m-1)}} \biggr)=f(x) \quad &\text{in}\; \Omega,\\ u=0 &\text{on}\; \partial \Omega,
 \end{cases}    
\end{equation}
with $\theta\ge0$ and $\Omega$ a bounded domain. 
In particular, they proved the existence and regularity of solutions, depending on summability of the datum $f$ reasoning by approximation to get a coercive differential operator on $W^{1,m}_0(\Omega)$.

We mention also \cite{BYM}, where Benkirane, Youssfi and Meskine prove existence and $L^\infty$-regularity for solutions of \eqref{pb_proto}. 

Following \cite{Ba}, we restrict our attention to (distributional) $C^1$ solutions of \eqref{probrad} satisfying $(i)_s$ and $(ii)_s$, under the following conditions on the parameters
\begin{enumerate}
\item[$(H_0)$]
$m>1$, $\alpha,\,\beta\in\mathbb{R}$ such that 
$$ N+\alpha-m>0, \qquad  \beta-\alpha+1>0,$$
\item[$(H_1)$]  
$m-1<p< m^*_{\alpha,\beta,\gamma}-1,\qquad \quad 0<\gamma<\Upsilon:=\dfrac{m(m-1)(\beta-\alpha+m)}{m(N+\beta+1)-N-\alpha},$
\end{enumerate}
where
$$m^*_{\alpha,\beta,\gamma}=\frac{m(m-1)(N+\beta)-\gamma[m(N+\beta+1)-N-\alpha]}{(m-1)(N+\alpha-m)}.$$
In particular,
$$m^*_{\alpha,\beta,\gamma}= m^*_{\alpha,\beta}-\frac{\gamma}{m-1}\bigl(m^*_{\alpha,\beta}-1)<m^*_{\alpha,\beta},$$ 
where $m^*_{\alpha,\beta}$ is given in \eqref{p*CDM},
since 
$\gamma<m-1$ from $\Upsilon<m-1$ by $(H_0)$ and being
\begin{equation}\label{first_sign}m(N+\beta+1)-N-\alpha>(m-1)(N+\alpha)>m(m-1)>0.
\end{equation}
We emphasize that $m^*_{\alpha,\beta,\gamma}$ has the role of the critical Sobolev exponent for problems with diffusion as \eqref{PROB}.
Indeed, it divides existence and nonexistence of positive radial solutions of problem \eqref{PROB}, as it will be clear from Corollary \ref{cor1} below.

Moreover, the range of $p$ in $(H_1)$ is nonempty since $\gamma$ stands below the upper threshold $\Upsilon$.
Furthermore, for $m=2$, the assumption $(H_1)$ becomes
$$1<p<2^*_{\alpha,\beta,\gamma}-1=2(1-\gamma)\frac{\beta-\alpha+2}{N+\alpha-2},
\qquad 0<\gamma<\frac{2(\beta-\alpha+2)}{N+2\beta-\alpha+2}, $$
we refer to \cite{Ba} for a similar assumption.
We point out that, assumptions $(H_0)$ and $(H_1)$, first appeared in \cite{CDM}, hold also in the case without weights, that is $\alpha=\beta=0$.

On the functions $a,g$ we assume
\begin{enumerate}
\item[$(H_2)$] $g:[0,\infty)\to[0,\infty)$ is a continuous nondecreasing function with $\displaystyle{\lim_{v\to\infty}\dfrac{g(v)}{v}=1}$.

\item[$(H_3)$]
$a:[0,\infty)\to(0,\infty)$ is a continuous function satisfying $c_1\le a(|x|)\le c_2$ in  $\mathbb R^N$, $c_1,$  $c_2$  positive constants.
\end{enumerate}

\smallskip We are ready to state the main theorems of our paper.
\begin{theorem}\label{th_1.1}
Assume $(H_0)$, $(H_1)$, $(H_2)$ and $(H_3)$. Then problem \eqref{PROB} has at least a positive radial solution.
\end{theorem}

This result extends and completes Theorem 1.1 in \cite{Ba} devoted to the case $m =2$. The technique used to prove Theorem \ref{th_1.1} is rather delicate and tangled because it requires several different tools. At the beginning, we need to investigate qualitative properties,  as well as the validity of a suitable  variation identity, for positive radial solutions to problem \eqref{liouv1} under 
$(H_1)'$ yielding a slightly different proof of the Liouville type result in \cite{CCM98, CM}. We emphasize
that $(H_1)'$ provides for $p$ a larger range than that in  $(H_1)$.

Secondly, under assumption $(H_1)$, we obtain nonexistence of positive radial solutions for the "broken problem"
\begin{equation}\label{liouv2}
\begin{cases}
      &-\mbox{div}(|x|^\alpha a_0(u)|\nabla u|^{m-2}\nabla u)=|x|^\beta u^p \quad\text{in} \quad \mathbb R^N\setminus\{0,\partial B_{s_0}\},\\
      & u(0)=1, \qquad  u(x)=1/d \,\, \text{ on}\quad\partial B_{s_0},\\ 
\end{cases}    
\end{equation}
for a certain  $d>0$ and $s_0>0$, where 
$$
a_0(u)=\begin{cases}
      1,\quad& \,\,|x|< s_0,\\
      1/(d^\gamma u^\gamma), & \,\,|x|\ge s_0.
\end{cases}
$$
Then, via blow up technique, adapted to our setting, we prove a priori uniform estimates for positive solutions of the parameterized truncated problem associated to \eqref{PROB}, for details see Section \ref{stima}, where in particular we prove that nonexistence occurs if the parameter is too large.
Next, a fixed point result by Krasnosel'skii  ensures the existence of  positive radial solutions for the truncated problem associated to \eqref{PROB}. We point out that these two results can be obtained under $(H_1)'$.  Finally, to obtain the existence of a positive radial solution for problem \eqref{PROB}, we apply a limit process to a  sequence of solutions of a suitable truncated problem together with an intensive  qualitative analysis and the use of Liouville type theorems.

The second part of the  paper is devoted to nonexistence results, aimed to give the counter part of Theorem \ref{th_1.1}, precisely we prove two nonexistence theorems, the first of them completes the picture together with Theorem \ref{th_1.1} and it is obtained through a variation identity for positive radial solutions of Poho\v zaev-Pucci-Serrin type.

For the nonexistence setting, we need to strengthen the regularity of $g$ and $a$, namely we assume
\begin{enumerate}
    \item [$(H_2)'$] $g:[0,\infty)\to[0,\infty)$ is a $C^1$ nondecreasing function.
    
    \item[$(H_3)'$] $a:[0,\infty)\to(0,\infty)$ is a $C^1$ nonincreasing   function.
    %function satisfying $$c_1\le a(|x|)\le c_2 \quad \text{for all}\quad x\in \mathbb R^N,$$where $c_1$ and $c_2$ are positive constants.
\end{enumerate}
In turn, we get 
\begin{theorem} \label{th_1.2} 
Assume $(H_0)$, $(H_2)'$ and $(H_3)'$.

Then,
problem \eqref{PROB} has no positive radial solutions for $p\ge m^*_{\alpha,\beta,\gamma}-1$ if either
\begin{enumerate}
\item[$(i)$] $\gamma<\Upsilon\,\,\quad\text{and}\quad m(N+\alpha-1)\ge \beta-\alpha+m\,\,$ or
\item[$(ii)$] $\gamma<\Upsilon_1\quad\text{and}\quad m(N+\alpha-1)< \beta-\alpha+m$,
\end{enumerate}
where 
\begin{equation}\label{ups1}
\Upsilon_1=\dfrac{m(m-1)\bigl[\alpha-\beta +m(N+\beta-1)\bigr](N+\alpha-1)}{[m(N+\beta+1)-N-\alpha][m(N+\alpha-1)-m+1]}
\end{equation}
and in case {\it (ii)} it results $\Upsilon_1<\Upsilon$.
\end{theorem}
The above theorem, new even in the case $m=2$, solves a problem left open by Theorem~1.2 in \cite{Ba} where  nonexistence  of positive radial solutions \eqref{PROB} with $m=2$  was obtained   and  $p\ge 2^*_{\alpha,\beta}-1$. Consequently 
in the interval 
$\bigl[2^*_{\alpha,\beta,\gamma}-1,\, 2^*_{\alpha,\beta}-1\bigr),$
the question of existence or nonexistence of solutions of problem \eqref{probrad} was still an open problem.

In  the prototype case $g(u)=u$ and $a(|x|)=1$, since assumptions $(H_2)'$ and $(H_3)'$ are trivially satisfied, combining Theorems \ref{th_1.1} and \ref{th_1.2} obtaining  
\begin{corollary}\label{cor1}
Assume $(H_0)$ and $\gamma<\Upsilon$. Then problem \begin{equation}\label{modello_caso}
\begin{cases}
  -\mbox{div}\biggl(\dfrac{|x|^\alpha|\nabla u|^{m-2}\nabla u}{\big(1+u\big)^\gamma} \biggr)=|x|^\beta u^p \quad &\text{in}\; B_R\setminus\{0\},\\ u=0 &\text{on}\; \partial B_R,
\end{cases}    
\end{equation}
has at least a positive radial solution if 
$ m-1<p<m^*_{\alpha,\beta,\gamma}-1$, while it admits no positive solutions if $p\ge m^*_{\alpha,\beta,\gamma}-1$ and 
$(i)$ or $(ii)$ in Theorem \ref{th_1.2}  hold.
\end{corollary}
In particular, in case $(i)$ of Theorem \ref{th_1.2}, it follows from Corollary \ref{cor1} that the critical value $ m^*_{\alpha,\beta,\gamma}$ divides existence from nonexistence, indeed  the following characterization holds.
\begin{corollary}
Assume $(H_0)$ with $p>m-1$ and $m(N+\alpha-1)> \beta-\alpha+m$. 

Then problem \eqref{modello_caso} has at least a positive radial solution if and only if $(H_1)$ holds.
\end{corollary}
Note that the case of same weights, that is $\alpha=\beta$, is covered by the above corollary when $m\ge2$, since condition $(i)$ in Theorem \ref{th_1.2} is equivalent to $N+\alpha-2\ge0$ which follows from $(H_0)$ if $m\ge 2$.

We observe that problem \eqref{modello_caso}, can be reduced to a problem with no diffusion, but with a different nonlinearity, indeed, for $\gamma\in(0,m-1)$, the following change of variable
$$w(x)=(1+u(x))^{1-\frac{\gamma}{m-1}}>1, $$
turns problem \eqref{modello_caso}  into
$$
\begin{cases}
  -\mbox{div}\bigl(|x|^\alpha|\nabla w|^{m-2}\nabla w\bigr)=\biggl(1-\dfrac{\gamma}{m-1}\biggr)^{m-1}|x|^\beta(w^{\frac{m-1}{m-1-\gamma}}-1)^p \quad &\text{in}\; B_R\setminus\{0\},\\ u=0 &\text{on}\; \partial B_R.
\end{cases}
$$

Finally, inspired by \cite{Ba}, we investigate also the case when $\beta-\alpha+m\le0$ in which roughly $m^*_{\alpha,\beta}$ loses its meaning becoming less that $m$, obtaining a second nonexistence result.
\begin{theorem}\label{th_1.3}
Assume $(H_3)$ and let $g$ be a continuous and nonnegative function. If
$$\gamma>0, \quad N+\beta>0 \quad \text{and}\quad \beta-\alpha+m\le0,$$
then problem \eqref{probrad} has no positive solutions.
\end{theorem}
This result extends Theorem 1.3 in \cite{Ba}.

\medskip The paper is organized as follows.
Section \ref{cap1} contains preliminary results, including regularity of positive radial solutions and classical tools used in the proofs of our main results, such as the fixed point Theorem by Krasnosel'skii. Then, in Section \ref{prep_lemmas}, by using a deep qualitative analysis in positive radial solutions, we prove  preparatory lemmas to the two Liouville type theorems given in Section \ref{liouv}. 
In Section \ref{stima}, we introduce and investigate the truncated and parameterized problem associated to \eqref{probrad} in order to obtain a priori estimates of positive solutions, a crucial tool in the proof of the main theorem. Section \ref{probtronc} is devoted to an existence theorem for the truncated problem obtained by a precise value of the parameter involved. Then, in Section \ref{exixtence} we prove the main existence result for positive radial solutions of problem \eqref{PROB} of the paper, Theorem \ref{th_1.1}. Finally, Section \ref{cap3} contains a new Poho\v aez-Pucci-Serrin type identity for positive radial solutions of \eqref{PROB} together with the proof of the two nonexistence results given by Theorems \ref{th_1.2} and \ref{th_1.3}.

\section{Classical results} \label{cap1}
In this section we recall some known results about regularity of positive radial solutions, together with qualitative 
properties, of problems of the type \eqref{PROB}. 
One of the pioneering papers in this direction is that of Ni and Serrin in \cite{NS}, where they  consider positive radial solutions of $\mbox{div}(\mathcal A(|\nabla u|)\nabla u)+f(u)=0$, namely solutions of 
\begin{equation}\label{nsprob}
    \begin{cases}
    -(r^{N-1}\mathcal A(|v'|)v')'=r^{N-1}f(v), \quad r>0,\\
    v(0)=v_0>0,
    \end{cases}
\end{equation}
for an operator of the form
\begin{equation}\label{operator_A}
\mathcal A:\mathbb R^+\to \mathbb R,\quad \mathcal A\in C^1(\mathbb R^+), \quad  t\mathcal A(t) \text{ strictly increasing},\quad
\lim_{t\to0^+}t\mathcal A(t)=0     
\end{equation}
and $f\in C(\mathbb R)$ with $f(0)=0$, $f>0$ in $\mathbb R^+$. Remarkable prototypes for $\mathcal A$ in $\mathbb R^+$ are the $m$-Laplacian operator $\mathcal A(t)=t^{m-2}$, $m>1$, 
the mean generalized curvature operator $\mathcal A(t)=(1+t^2)^{m/2 -1}$, $m\in (1,2]$, and the mean curvature operator $A(t)=(1+t^2)^{-1/2}$.

In particular,  they prove in Proposition 1 in \cite{NS} that every solution $v=v(r)$ of \eqref{nsprob}
is continuously differentiable on some interval $0\le r\le R$ with $v'(0)=0$. 
The proof of the above result is based on the application of Schauder's fixed point theorem, applied to the compact operator
\begin{equation}\label{T_operator} T[v](r)=v_0-\int_0^r\varphi\biggl(\int_0^tf(v(s))\biggl(\frac st\biggr)^{N-1}ds\biggr)dt,
    \end{equation}
with $r\ge0$ and $v\in C[0,R]$, for some $R>0$, where $\varphi$ is the inverse function of $t\mathcal A(t)$, with $\varphi(0)=0$, in 
 $\mathfrak C=\bigl\{v\in C([0,R]): \|v(r)-v_0\|_\infty\le v_0/2\bigr\},$ for $R$ suitably small so that the value
$\int_0^tf(v(s))( s/t)^{N-1}ds$ small, for all  $t\in (0,R]$,
gives $T(\mathfrak C)\subset \mathfrak C$.

Furthermore, they proved that 
positive solutions of \eqref{nsprob}  are of class $C^2$, namely they are classical, as far as $v'(r)\ne0$.  

Moving to radial problems with different weights, problem \eqref{nsprob} changes as follows
\begin{equation}\label{pesi_diversi}
 \begin{cases}
  -\left(r^{N+\alpha-1}\mathcal A(|v'|)v'\right)'=r^{N+\beta-1}f(v(r)), \quad r>0,
  \\v(0)=v_0>0,
   \end{cases}
\end{equation}
thus the operator $T$ defined in \eqref{T_operator} needs to be replaced by
$$T[v](r)=v_0-\int_0^r\varphi\biggl(
\frac1 {t^{N+\alpha-1}}\int_0^t f(v(s)) s^{N+\beta-1}ds\biggr)dt,$$
so that $T(\mathfrak C)\subset \mathfrak C$ holds if
$$0\le \frac1 {t^{N+\alpha-1}}\int_0^t f(v(s)) s^{N+\beta-1}ds\le \frac{t^{\beta-\alpha+1}}{N+\beta}\max_{t\in[0,R]}f(v(t))$$
is sufficiently small for all $t\in[0,R]$, with $R>0$ suitable. Thus, as observed in \cite{CDM}, condition $\beta-\alpha+1>0$ is necessary to obtain regularity $C^1[0,R]$ of positive solutions of \eqref{pesi_diversi}, as well as $v'(0)=0$.

Furthermore, for any nonnegative solution $v$ of \eqref{pesi_diversi} we have 
\begin{equation}\label{v'_eq}
-\mathcal A(|v'(r)|)v'(r)=\frac{1}{r^{N+\alpha-1}}\int_0^r s^{N+\beta-1}f(v(s))ds, \quad r\in (0, R]
\end{equation}
and from the positivity of the right hand side, we deduce that  $v'(r)<0$ for all $0<r\le R$, being  $\mathcal A>0$ in $\mathbb R^+$, so that  regularity $C^2(0,R)$ of $v$ follows using also $C^1$ regularity of $\mathcal A$.

We point out that in the model case $f(v)=v^p$, $v>0$, thus, assuming $\beta-\alpha+1>0$,  then  $v$ can be continued for $r>R$ by the boundedness of $v$ and $v'$. Indeed, $v$ is positive and bounded by  monotonicity, while  the boundedness of $v'$ is a consequence of 
$$\frac{1}{r^{N+\alpha-1}}\int_0^r s^{N+\beta-1}v(s)^pds\le \frac{v_0^p}{N+\beta}r^{\beta-\alpha+1}$$
so that from \eqref{v'_eq}, since $v'(r)<0$,  and the  strictly increasing monotonicity of $\varphi$,
we arrive to 
$$|v'(r)|\le \varphi \biggl(\frac{v_0^p}{N+\beta}r^{\beta-\alpha+1}\biggr).$$
In turn $v\in C^1[0,\infty)$ and by  $v'<0$ we get $v\in C^2(\mathbb R^+;\mathbb R^+)$.
For a more complete existence result we refer to Theorem $5.2$ in \cite{CDM}.

Now, returning to problem \eqref{probrad} and arguing as above, we discuss the regularity of positive solutions.
Following \cite{GMPS}, we give the definition of weak (distribution) solutions of \eqref{PROB},
that is equivalent to the notion of semiclassical  and classical $C^1$ solution, since we consider positive solutions of problem \eqref{PROB} with continuous nonlinearities which vanishes at $0$, for details see Proposition 2.1 in \cite{GMPS}.

\begin{definition}
A weak (distribution) solution of \eqref{PROB} is a nonnegative function $u$ of $C^1(B_R)\cap C({\overline{B_R}})$  which verifies
$$
\int_{B_R}\frac{|x|^\alpha|\nabla u|^{m-2}\nabla u}{(a(|x|+g(u))^\gamma}\cdot\nabla\psi dx=\int_{B_R}|x|^\beta u^p\psi dx
$$
for all $C^1$ functions $\psi=\psi(x)$ with compact support in $B_R$.
\end{definition}

Equivalently for the radial case, $v$ is a weak (distribution) solution of \eqref{probrad} if  $v$ is a nonnegative function $v$ of $C^1(0,R)\cap C[0,R]$  satisfying
$$
\int_0^R\frac{r^{N+\alpha-1}|v'|^{m-2}v'}{(a(r)+g(v))^\gamma}\cdot\psi'dr=\int_0^Rr^{N+\beta-1}v^p(r)\psi dr
$$
for all $C^1$ functions $\psi=\psi(r)$ with compact support in $[0,R)$. 
Then, by  distribution arguments, we have that $v$ satisfies
$$-\frac{r^{N+\alpha-1}|v'(r)|^{m-2}v'(r)}{(a(r)+g(v(r)))^\gamma}=\int_0^rs^{N+\beta-1}v^p(s)ds
$$
in $[0,R)$. Since the right hand side is continuously differentiable in $r$, it follows that
$r^{N+\alpha-1}|v'|^{m-2}v'/(a(r)+g(v))^\gamma\in C^1(0,R)$
so that $v$ is a classical solution of problem \eqref{PROB}.
Thus, in our case the definition of weak solution is compatible with that of classical solution, indeed if  $v\in C^1(0,R)\cap C[0,R]$ is a weak positive solution of \eqref{probrad}, then
$$
\frac{(a(r)+g(v(r)))^\gamma}{r^{N+\alpha-1}}\int_0^rs^{N+\beta-1}v^p(s)ds>0, \quad \text{in}\,\, (0,R),
$$
from $a>0$, $g\ge0$ and $v$ positive in $(0,R)$, so that $-|v'(r)|=v'(r)<0$ in $(0,R)$ and in turn
\begin{equation}\label{v'_prob}
v'(r)=-\biggl(\frac{(a(r)+g(v(r)))^\gamma}{r^{N+\alpha-1}}\int_0^rs^{N+\beta-1}v^p(s)ds\biggl)^{1/(m-1)} \quad \text{in} \,\, (0,R).    
\end{equation}
In addition, $v$ is $C^1$ at $r=0$, so that  \eqref{v'_prob} holds in $[0,R]$. Indeed, by $g$ nondecreasing and by $(H_3)$, then
$c_1+g(0)\le a(r)+g(v(r))\le c_2+g(v(0))$, consequently
using that $\beta-\alpha+1>0$ and L'H\^opital's rule, we arrive to
$$
\lim_{r\to0^+}|v'(r)|\le  C \lim_{r\to 0^+} r^{\frac{\beta-\alpha+1}{m-1}}=0,\qquad C>0.
$$
In conclusion,  $v\in C^1[0,R]$ with 
\begin{equation}\label{derivatasegno}
v'(0)=0 \quad\text{and}\quad v'(r)<0\quad\text{for all}\,\, 0<r\le R. 
\end{equation}
In turn,  $v\in C^2(0,R)$  being $\varphi(t)=t^{1/(m-1)}\in C^1(0,\infty)$ in \eqref{v'_prob}. 
Summarizing, under the assumptions $(H_3)$,  $g$ continuous and nondecreasing and $\beta-\alpha+1>0$,
we obtain that every positive weak solution $v$ of \eqref{probrad} is 
\begin{equation}\label{reg2}
v\in C^1[0,R]\cap C^2(0,R).    
\end{equation}

We end this section by stating a fixed point Theorem of Krasnosel'skii in \cite{K}, aimed to obtain existence of solutions of problem \eqref{PROB} as fixed points of compact operators defined in a cone, which is one of the crucial tools in the proof of  Theorem \ref{th_1.1}. 

\begin{theorem}[Krasnosel'skii] \label{kran_th}
Let $\mathcal{K}$ be a cone in a Banach space, and let $F:\mathcal{K}\to\mathcal{K}$ be a compact operator such that $F(0)=0$. Suppose there exists $\delta>0$ verifying
\begin{enumerate}
    \item [$(a)$] \quad $u\ne tF(u)$, for all $\|u\|=\delta$ and $t\in[0,1]$.
\end{enumerate}    
Suppose further that there is a compact homotopy $\mathcal H:[0,1]\times\mathcal{K}\to\mathcal{K}$ and $\eta>\delta$ such that:
\begin{enumerate}
    \item[$(b)$] \quad $F(u)=\mathcal H(0,u)$, for all $u\in\mathcal{C}$. 
    \item[$(c)$] \quad $\mathcal H(t,u)\ne u$, for all $\|u\|=\eta$ and $t\in[0,1]$.
    \item[$(d)$] \quad $\mathcal H(1,u)\ne u$, for all $\|u\|\le\eta$.
\end{enumerate}
Then $F$ has a fixed point $u_0$ verifying $\delta<\|u_0\|<\eta$.
\end{theorem}

\section{Preparatory Lemmas} \label{prep_lemmas}

In this section we prove two preparatory lemmas dealing with the following related problem 
\begin{equation} \label{p1}
   \begin{cases}
   -(r^{N+\alpha-1}|u'(r)|^{m-2}u'(r))'=\lambda r^{N+\beta-1}u^{\wp}(r),& r\in (s_{0},\infty), \\u(s_{0})=1, \qquad  u'(s_{0})\leq 0.
   \end{cases}
\end{equation}
with $s_{0}\geq 0$, $\lambda >0$ and $\wp$ verifying 
$(H_1)'$.

\begin{lemma} \label{lemma2.1}
Let $u\in C^1(s_0,\infty)$ be a nonnegative solution of \eqref{p1}, with $s_{0}\geq 0$, where $\alpha, \beta, m$ satisfy $(H_{0})$ and take $\wp>0$. Then the function $U_\varrho(r):=ru'(r)+\varrho u(r)$ is nonnegative and nonincreasing for 
\begin{equation}\label{def_ro}
\varrho=\frac{N+\alpha-m}{m-1}.
\end{equation}
In particular $r^\varrho u(r)$ is nondecreasing on $(s_{0},\infty)$.
\end{lemma}
\begin{proof}
From the equation in \eqref{p1} we observe that 
\begin{equation}\label{udec}
u'(r)< 0\quad\text{in}\quad  (s_{0},\infty)
\end{equation}
for every $u$ nonnegative solution of \eqref{p1}. Indeed, integrating the equation from $s_{0}$ to $r$ with $r>s_{0}$, we obtain
$$
 -r^{N+\alpha-1}|u'(r)|^{m-2}u'(r)+s_{0}^{N+\alpha-1}|u'(s_{0})|^{m-2}u'(s_{0})=\lambda \int_{s_{0}}^{r}{s^{N+\beta-1}u^{\wp}(s)ds}> 0,$$ 
so we have 
$$-r^{N+\alpha-1}|u'(r)|^{m-2}u'(r)>-s_{0}^{N+\alpha-1}|u'(s_{0})|^{m-2}u'(s_{0})\geq 0\quad\text{in}\quad  (s_{0},\infty)$$
being $u'(s_{0})\leq0$ since $u$ is a solution of \eqref{p1}, so that \eqref{udec} is proved and $u\in C^2(s_0,\infty)$ arguing as in Section \ref{cap1}, cfr. \eqref{reg2}.

From the nonnegativity of $u$, we know that $-(r^{N+\alpha-1}|u'(r)|^{m-2}u'(r))'\geq 0$ for all $r\in (s_{0},\infty)$, namely, 
by \eqref{udec},
\begin{align*}
    (N+\alpha-1)r^{N+\alpha-2}|u'(r)|^{m-1}-r^{N+\alpha-1}(m-1)|u'(r)|^{m-2}u''(r)\geq 0,    
\end{align*}
for all $r\in(s_0,\infty)$. Consequently, multiplying the last inequality by $u'(<0)$,  we obtain
$$r^{N+\alpha-2}\left[\dfrac{N+\alpha-1}{m-1} u'(r)+r u''(r)\right]\leq0 \quad\text{in}\quad  (s_{0},\infty),
$$
being $r>s_0\geq0$ and $m>1$, so that, by \eqref{def_ro}, we get
$$
[U_\varrho(r)]'
    =(1+\varrho)u'+ru''(r)=\dfrac{N+\alpha-1}{m-1} u'(r)+r u''(r)\leq0$$
for all $r\in (s_{0},\infty)$. Thus, $U_\varrho(r)$ is a nonincreasing function on $(s_0,\infty)$.

On the other hand, to obtain the nonnegativity of $U_\varrho(r)$, we argue by contradiction by assuming   that there exists $r_0>s_{0}$ and $m_{0}<0$ such that $U_\varrho(r_0)\leq m_{0}$. By the fact that $U_\varrho(r)$ is a nonincreasing function on $(s_{0},\infty)$, 
then $U_\varrho(r)\leq m_{0}$ for all $r\geq r_0$, that is
$$ru'(r)\leq m_{0}-\varrho u(r)\leq m_{0},$$
since $\varrho>0$ by $(H_0)$ and the solution $u$ is nonnegative. Thus, for all $r\in(s_{0},\infty)$ we get
$$u'(r)\leq m_{0}r^{-1}$$
and integrating from $r_0$ to $r$, we have that 
$$ u(r)-u(r_0)\leq m_{0}\ln{\dfrac{r}{r_0}},$$
yielding
$\lim_{r \to \infty}u(r)=-\infty,$
since $m_0<0$. Thus the required contradiction is reached since $u$ is a nonnegative function, yielding that  $U_\varrho\ge 0$  in $(s_0,\infty)$. In particular, it follows that for all $r\in(s_{0},\infty)$
\begin{align*}
    \bigl(r^{\varrho}u(r)\bigr)'&=r^{\varrho}u'(r)+\varrho r^{\varrho-1}u(r)=r^{\varrho-1}U_\varrho (r) \geq 0.
\end{align*}
The proof is so concluded.
\end{proof}

In the next lemma we will show a priori estimates for solutions of \eqref{p1} and a variational identity for positive radial solutions of problem \eqref{p1}, both  crucial ingredients to obtain Liouville type results in the next section.

\begin{lemma} \label{lemma2.2}
Assume $(H_{0})$, $\wp>m-1$ and let $u$ be a positive solution of \eqref{p1}. 

Then, there are two positive constants $C_i=C(N,\alpha,\beta,\lambda,\wp, m)$ for $i=1,2$, such that
\begin{equation}\label{tesi1_Lemma2.2}
r^{N+\beta}u^{\wp+1}(r)\leq \begin{cases}
    C_1r^{-\frac{m(N+\beta)-(N+\alpha-m)(\wp+1)}{\wp-m+1}} \qquad &\text{if}\,\,N+\beta\neq\varrho\wp,\\
    C_2r^{-\frac{(N+\alpha-m)(\wp+1)-(m-1)(N+\beta)}{m-1}}&\text{if}\,\,N+\beta=\varrho\wp,
\end{cases}
\end{equation}
in $(s_0,\infty)$ with $s_0\ge0$ and where $\varrho$ is defined in \eqref{def_ro}. 

Moreover, the following identity holds
\begin{equation}\begin{aligned}\label{tesi2_Lemma2.2}
    &\lambda\left(\dfrac{N+\beta}{\wp+1}-\dfrac{N+\alpha-m}{m}\right)\int_{s_0}^{r}{s^{N+\beta-1}u^{\wp+1}(s)ds}\\
    &=-\frac{m-1}m r^{N+\alpha-1}|u'(r)|^{m-1}
    U_\varrho(r)
    +\frac{m-1}m s_0^{N+\alpha-1}|u'(s_0)|^{m-1}U_\varrho(s_0)\\&\quad+
    \dfrac{\lambda}{\wp+1}[r^{N+\beta}u^{\wp+1}(r)-s_0^{N+\beta}],
\end{aligned}\end{equation}
 in $(s_0,\infty)$, where $U_\varrho(r)$ is defined in Lemma $\ref{lemma2.1}$.
\end{lemma}
\begin{proof}
Assume that $u$ is a  positive solution of  problem \eqref{p1}. Fix $r, t$ such that $s_0< r<t$, then integrate the equation in \eqref{p1} from $r$ to $t$, and use 
 $r^{N+\alpha-1}|u'(r)|^{m-1}\ge0$, so that
$$t^{N+\alpha-1}|u'(t)|^{m-1}\ge\lambda\int_{r}^{t}{s^{N+\beta-1}u^{\wp}(s)ds}$$
thus, by the monotonicity of $s^\varrho u(s)$ given by Lemma \ref{lemma2.1} and \eqref{udec}, we observe that the following holds
\begin{equation}
\begin{aligned}\label{ineq1}
&t^{N+\alpha-1}|u'(t)|^{m-1}\geq\lambda\int_{r}^{t}{s^{N+\beta-1-\varrho \wp}(s^\varrho u(s))^{\wp} ds}
\geq\lambda(r^\varrho u(r))^{\wp}\!\int_{r}^{t}{s^{N+\beta-1-\varrho \wp}ds}\\&=\lambda r^{N+\beta}u^{\wp}(r)\begin{cases}\dfrac1{N+\beta-\varrho \wp}\biggl[\biggl(\dfrac t r\biggr)^{N+\beta-\varrho \wp}-1\biggr] \qquad &\text{if}\,\, N+\beta\neq \varrho \wp,\\\\
r^{\varrho\wp-N-\beta}\log\dfrac tr \qquad &\text{if}\,\, N+\beta= \varrho \wp.\end{cases}
\end{aligned}\end{equation}
%we note that the right hand side is always positive.
We consider separately the two cases starting with $N+\beta\neq\varrho\wp$.
Thus, observing that the right hand side of \eqref{ineq1} is always positive,  we have
\begin{equation}\label{magg_1}
\begin{aligned}
t^{(N+\alpha-1)/(m-1)}|u'(t)|\ge\biggl(\lambda r^{\varrho p}u^{\wp}(r)\dfrac{t^{N+\beta-\varrho \wp}-r^{N+\beta-\varrho \wp}}{N+\beta-\varrho \wp}\biggr)^{\frac{1}{m-1}}.
\end{aligned}
\end{equation}
Now using \eqref{udec} and the nonnegativity of $U_\varrho$ given by Lemma \ref{lemma2.1}, we get
\begin{equation} \label{magg_2}
t^{(N+\alpha-1)/(m-1)}|u'(t)|=t^{(N+\alpha-m)/(m-1)}t|u'(t)|\le t^{(N+\alpha-m)/(m-1)} \varrho u(t),
\end{equation}
so that, combining \eqref{magg_1} and \eqref{magg_2} we arrive to 
$$t^{(N+\alpha-m)/(m-1)}\varrho u(t)\geq\left(\lambda r^{\varrho \wp}\frac{t^{N+\beta-\varrho \wp}-r^{N+\beta-\varrho \wp}}{N+\beta-\varrho \wp}u^{\wp}(r)\right)^{\frac{1}{m-1}},\quad t>r.$$
Taking $t=2r$ and using the fact that $u$ is a positive and decreasing function, it holds
$$
\begin{aligned}
(2r)^{(N+\alpha-m)/(m-1)}\varrho u(r)&\geq
    (2r)^{(N+\alpha-m)/(m-1)}\varrho u(2r)\\
    &\geq\left[\dfrac{\lambda(2^{N+\beta-\varrho \wp}-1)}{N+\beta-\varrho \wp}r^{N+\beta}\right]^{\frac{1}{m-1}}u^{\frac{\wp}{m-1}}(r)
\end{aligned}
$$
from which we have
$$u^{\frac{\wp-m+1}{m-1}}(r)\leq Cr^{\frac{\alpha-\beta-m}{m-1}},  $$
where $\alpha-\beta-m<\alpha-\beta-1<0$ by ($H_0$).
In particular, elevating both members to $(m-1)(\wp+1)/(\wp-m+1)>0$, we have, for $N+\beta\neq\varrho\wp$,
\begin{equation}\label{maggior_u}
u^{\wp+1}(r)\leq C r^{-\frac{(\beta-\alpha+m)(\wp+1)}{\wp-m+1}}, \quad r\in(s_0,\infty).    
\end{equation}
Now we consider the case $N+\beta=\varrho\wp$ and we argue as above. 
By \eqref{ineq1} we get \begin{equation}\label{magg_1bis}
\begin{aligned}
t^{(N+\alpha-m)/(m-1)}\varrho u(t)\geq t^{(N+\alpha-1)/(m-1)}|u'(t)|\ge\biggl(\lambda r^{\varrho p}u^{\wp}(r)\log\dfrac{t}{r}\biggr)^{\frac{1}{m-1}},\quad t>r,
\end{aligned}
\end{equation}
where we have used \eqref{magg_2} which holds also in this case.
Taking again $t=2r$, we obtain
$$
\begin{aligned}
(2r)^{(N+\alpha-m)/(m-1)}\varrho u(r)&\geq
    (2r)^{(N+\alpha-m)/(m-1)}\varrho u(2r)\geq\left[\lambda\log2r^{\varrho \wp}\right]^{\frac{1}{m-1}}u^{\frac{\wp}{m-1}}(r)
\end{aligned}
$$
yielding
$$u^{\frac{\wp-m+1}{m-1}}(r)\leq Cr^{\frac{N+\alpha-m}{m-1}\bigl(1-\frac{\wp}{m-1}\bigr)},$$
where the exponent of $r$ is negative by $(H_0)$ and $\wp>m-1>0$.
In particular, elevating both members to $(m-1)(\wp+1)/(\wp-m+1)>0$, we have
\begin{equation}\label{maggior_ubis}
u^{\wp+1}(r)\leq C r^{-\frac{(N+\alpha-m)(\wp+1)}{m-1}}, \quad r\in(s_0,\infty).    
\end{equation}
Thus, \eqref{tesi1_Lemma2.2} follows multiplying \eqref{maggior_u} and \eqref{maggior_ubis} by $r^{N+\beta}$.

Now we are ready to prove identity \eqref{tesi2_Lemma2.2}. First, multiply equation in \eqref{p1} by $ru'(r)$ and then integrate from $s_0$ to $r$ we get
\begin{equation}\label{int_1}
    -\int_{s_0}^{r}su'(s)\left(s^{N+\alpha-1}|u'(s)|^{m-2}u'(s)\right)'ds=\int_{s_0}^{r}{\lambda s^{N+\beta}u^{\wp}(s)u'(s)ds}.
\end{equation}
We analyze separately the left hand side  and the right hand side  of \eqref{int_1}. We start by integrating  by parts the left hand side 
\begin{equation} \label{primo_membro}
\begin{aligned}
   \int_{s_0}^{r}su'(s)&\left(s^{N+\alpha-1}|u'(s)|^{m-2}u'(s)\right)'ds
   =\frac{m-1}mr^{N+\alpha}|u'(r)|^{m}\\&-\frac{m-1}ms_0^{N+\alpha}|u'(s_0)|^m
   +\frac{N+\alpha-m}m\int_{s_0}^{r}s^{N+\alpha-1}|u'(s)|^{m}ds.
\end{aligned}    
\end{equation}
Furthermore, integrating by parts the right hand side of \eqref{int_1} and using $u(s_0)=1$, we have
\begin{equation} \label{secondo_membro}
\begin{aligned}
\int_{s_0}^{r}{\lambda s^{N+\beta}u^{\wp}(s)u'(s)ds}=&\dfrac{\lambda}{\wp+1}
\left[r^{N+\beta}u^{\wp+1}(r)-s_0^{N+\beta}\right]\\
&\quad-\frac{\lambda(N+\beta)}{\wp+1}\int_{s_0}^{r}{s^{N+\beta-1}u^{\wp+1}(s)ds},
\end{aligned}
\end{equation}
so that, by \eqref{primo_membro} and \eqref{secondo_membro}, the equality \eqref{int_1} becomes
\begin{equation} \label{int_1_final}
\begin{aligned}
\frac{m-1}m&\biggl[s_0^{N+\alpha}|u'(s_0)|^m-
r^{N+\alpha}|u'(r)|^{m}\biggr]-\frac{N+\alpha-m}m\int_{s_0}^{r}{s^{N+\alpha-1}|u'(s)|^{m}ds}\\
&=\dfrac{\lambda}{\wp+1}\left[r^{N+\beta}u^{\wp+1}(r)-s_0^{N+\beta}-(N+\beta)\int_{s_0}^{r}{s^{N+\beta-1}u^{\wp+1}(s)ds}\right].
\end{aligned}\end{equation}
On the other hand, multiplying the equation in \eqref{p1} by $u$, then integrating it from $s_0$ to $r$ and using again that $u(s_0)=1$, we obtain
\begin{equation} \label{integral}
    \begin{aligned}
    \int_{s_0}^{r}&{s^{N+\alpha-1}|u'(s)|^{m}ds}=-r^{N+\alpha-1}|u'(r)|^{m-1}u(r)\\
    &\qquad\qquad+s_0^{N+\alpha-1}|u'(s_0)|^{m-1}
    +\lambda\int_{s_0}^{r}{s^{N+\beta-1}u^{\wp+1}(s)ds}.
    \end{aligned}
\end{equation}
Replacing \eqref{integral} in \eqref{int_1_final}  we obtain
\begin{align*}
    &-\frac{m-1}m r^{N+\alpha}|u'(r)|^{m}+\dfrac{m-1}m s_0^{N+\alpha}|u'(s_0)|^{m}\\
    &-\frac{N+\alpha-m}{m}\left[-r^{N+\alpha-1}|u'(r)|^{m-1}u(r)+s_0^{N+\alpha-1}|u'(s_0)|^{m-1}+\lambda\int_{s_0}^{r}{s^{N+\beta-1}u^{\wp+1}(s)ds}\right]\\
    &\qquad=\dfrac{\lambda}{\wp+1}\biggl[r^{N+\beta}u^{\wp+1}(r)-s_0^{N+\beta}-(N+\beta)\int_{s_0}^{r}{s^{N+\beta-1}u^{\wp+1}(s)ds}\biggr].
\end{align*}
Hence, we have
$$\begin{aligned} 
\lambda&\left(\dfrac{N+\beta}{\wp+1}-\dfrac{N+\alpha-m}{m}\right)\int_{s_0}^{r}{s^{N+\beta-1}u^{\wp+1}(s)ds}\\&
\quad=\frac{m-1}m r^{N+\alpha}|u'(r)|^{m}
-\frac{N+\alpha-m}{m}r^{N+\alpha-1}|u'(r)|^{m-1}u(r) +\frac{\lambda}{\wp+1}r^{N+\beta}u^{\wp+1}(r)\\&
\qquad
-\frac{m-1}m s_0^{N+\alpha}|u'(s_0)|^{m}+\frac{N+\alpha-m}{m}s_0^{N+\alpha-1}|u'(s_0)|^{m-1}
-\dfrac{\lambda}{p+1} s_0^{N+\beta},
\end{aligned}$$
so that \eqref{tesi2_Lemma2.2} is proved.
\end{proof}

Consequently, the following property holds for positive solutions of \eqref{p1}.

\begin{corollary} \label{cor_p}
Assume $(H_{0})$ and $(H_{1})'$. Then, every  positive solution of \eqref{p1}
is such that
\begin{equation} \label{infinitesimo}
    \lim_{r\to\infty}r^{N+\beta}u^{\wp+1}(r)=0.
\end{equation} 
In particular, $ \lim_{r\to\infty}u(r)=0$.
\begin{proof}
First of all, we observe that to obtain \eqref{infinitesimo}, we need to verify that the exponents of $r$ in \eqref{tesi1_Lemma2.2} are negative, namely $$m(N+\beta)-(N+\alpha-m)(\wp+1)>0, \quad \wp-m+1>0\qquad\text{if}\quad N+\beta\neq \varrho\wp$$
and 
$$(N+\alpha-m)(\wp+1)-(m-1)(N+\beta)>0\qquad\text{if}\quad N+\beta=\varrho\wp.
$$
Actually, the first case follows from $(H_1)'$ using that  $\wp+1<m^*_{\alpha,\beta}$,
while the second, which can occur since 
$\beta-\alpha+m>0$ gives
$$
\wp=\frac{N+\beta}{\varrho}=\frac{(N+\beta)(m-1)}{N+\alpha-m}\in\biggl(m-1,\frac{(m-1)(N+\beta)+\beta-\alpha+m}{N+\alpha-m}\biggr)= \bigl(m-1,m^*_{\alpha,\beta}-1\bigr),
$$ 
can be verified directly
being, by $(H_0)$,
$$
(N+\alpha-m)(\wp+1)-(m-1)(N+\beta)=
N+\alpha-m>0$$
as soon as  $N+\beta=\varrho\wp$. 

In turn, by letting $r\to\infty$ in the a priori estimate  \eqref{tesi1_Lemma2.2}, in both cases property 
\eqref{infinitesimo} immediately follows.
\end{proof}
\end{corollary}

\section{Liouville type theorems} \label{liouv}

In this section we prove two  Liouville type theorems for equations of type \eqref{PROB} in $\mathbb R^N$ with $\gamma=0$ and of type \eqref{liouv2}, respectively.  
 Actually, the first of them,  which deals with the subcase of problem \eqref{p1} when $u'(s_0)=0$ and $s_0=0$,  is an application of Theorem 7.3 in \cite{CM} with $a(r)=\lambda$, which extends to the case $1<m<2$ Theorem 3.2 in \cite{CCM98}, where nonexistence was obtained only for $m\ge2$, for details see Remark \ref{cmm} below.

For completeness,  in what follows we give the proof,  inspired to \cite{Ba} for the case $m=2$, of the first Liouville result, since it is slightly different 
respect to \cite{CM}, indeed it based on the use of   Lemma \ref{lemma2.1} and \ref{lemma2.2} given in Section \ref{prep_lemmas}.

\begin{theorem} \label{teorema1}
Assume $(H_0)$, $(H_1)'$ and let $\lambda>0$. Then, the following problem 
\begin{equation}
\begin{cases} \label{p2}
   -(r^{N+\alpha-1}|u'(r)|^{m-2}u'(r))'=\lambda r^{N+\beta-1}u^{\wp}(r),& r\in (0,\infty), \\u(0)=1, \qquad  u'(0)=0,
   \end{cases}
\end{equation}
does not admit  positive solutions.
\end{theorem}
\begin{proof}
We argue by contradiction by assuming   that there exists positive solution $u$ of the problem \eqref{p2}. Then, we apply Corollary \ref{cor_p} Lemma \ref{lemma2.2}  for $u$ with $s_0=0$, so that \eqref{infinitesimo} holds and 
using that $u'(0)=0$, the identity \eqref{tesi2_Lemma2.2}
becomes
\begin{align*}
    &\lambda\left(\dfrac{N+\beta}{\wp+1}-\dfrac{N+\alpha-m}{m}\right)\int_{0}^{r}{s^{N+\beta-1}u^{\wp+1}(s) ds}\\
    &=-\frac{m-1}m r^{N+\alpha-1}|u'(r)|^{m-1} U_\varrho(r)+
   \dfrac{ \lambda}{\wp+1}r^{N+\beta}u^{\wp+1}(r)
\end{align*}
for any $r\in(0,\infty)$.
Notice that Lemma \ref{lemma2.1} holds, thus we know that 
$$ r^{N+\alpha-1}|u'(r)|^{m-1} U_\varrho(r)\ge0,$$
so we have
$$\left(\dfrac{N+\beta}{\wp+1}-\dfrac{N+\alpha-m}{m}\right)\int_{0}^{r}{s^{N+\beta-1}u^{\wp+1}(s) ds}\le \dfrac{1}{\wp+1}r^{N+\beta}u^{\wp+1}(r)$$
for any $r\in(0,\infty)$. From the last inequality, letting $r\to\infty$ and using \eqref{infinitesimo}, we obtain
\begin{equation} \label{disuguaglianza}
\left(\dfrac{N+\beta}{\wp+1}-\dfrac{N+\alpha-m}{m}\right)
\int_{0}^{\infty}{s^{N+\beta-1}u^{\wp+1}(s) ds}\le0.   
\end{equation}
In particular, by $(H_1)'$ we have that the coefficient of the integral in \eqref{disuguaglianza} is positive, so that
$$\int_{0}^{\infty}{s^{N+\beta-1}u^{\wp+1}(s) ds}\le0,$$
which is the required contradiction, since $u$ is a positive solution of \eqref{p2}.
\end{proof}

\begin{remark}\label{cmm}
We point out that the proof technique both of Theorem 7.3 in \cite{CM} and of Theorem 3.2 in \cite{CCM98} is different from the one above. Indeed, in both cases, the most delicate part in their proof consists in proving the following estimate
\begin{equation}\label{fincc}
u(r)\le C r^{-(N+\alpha-m)/m}\quad \text{in}\,\, [0,\infty),
\end{equation}
which was obtained in \cite{CCM98} only for $m\ge2$ and then extended to $1<m<2$ by Caristi and Mitidieri in \cite{CM}. Actually \eqref{fincc} is reached for $m\ge2$ in 
\cite[Proposition 2.2, Lemma 2.1]{CCM98} by proving the nonnegativity of a certain integral function, in our case their function $a(\xi)$ is a positive constant,   given by 
\begin{equation}\label{liminf}
\liminf_{r\to\infty} r^{(N+\alpha-m)(\wp+1)/(m-1)m} \int_r^\infty \xi^{\frac{(m-1)(N+\beta)-N-\alpha+1-\wp(N+\alpha-m)}{m-1}}d\xi\ge0.
\end{equation}
While Caristi and Mitidieri,  
in order to extend \eqref{fincc} to the case $1<m<2$,  use \cite[Lemma 7.4]{CM} and the nonnegativity of the following function
$$\liminf_{r\to\infty} r^{(N+\alpha-m)(\wp-m+1)/(m-1)m} \int_r^\infty \xi^{N+\beta-1-\frac{N+\alpha-m}{m-1}\wp}d\xi\ge0.$$
\end{remark}

\medskip
Now we prove a second nonexistence result for the "broken problem" defined as follows
\begin{equation} \label{P}
    \begin{cases}
   \begin{cases}
      -(r^{N+\alpha-1}|u'(r)|^{m-2}u'(r))'= r^{N+\beta-1}u^{p}(r),& r\in (0,s_0), \\u(0)=1, \qquad u'(0)=0,
   \end{cases}\\
   \begin{cases}
      -\left(\dfrac{r^{N+\alpha-1}|u'(r)|^{m-2}u'(r)}{d^\gamma u^\gamma(r)}\right)'= r^{N+\beta-1}u^{p}(r),&\quad r\in (s_0,\infty), \\u(s_0)=d^{-1},
   \end{cases}
   \end{cases} 
\end{equation}
where $d>1$ and $s_0>0$. In particular, we note that the equation in the second system coincides with the equation in the first one as $r\to s_0^+$ by the definition of $u(s_0)$. Moreover, the solution restricted to the first interval is $C^1[0,s_0]$ as we have seen in Section \ref{cap1} and in addition  $u'(s_0)<0$, so that   this latter condition is tacitly assumed in the second system.

We emphasize, that in the next nonexistence result we need to strengthen
 assumption $(H_1)'$, by assuming condition $(H_1)$. Precisely, 
\begin{theorem} \label{teorema2}
Assume $(H_0)$ and $(H_1)$. Then, problem \eqref{P} does not admit any positive solution.
\end{theorem}

\begin{proof}
We argue by contradiction, by assuming that there exists a positive solution $u$ of the problem \eqref{P}. Considering the following change of variables
$$v(r)=\bigl( d\cdot u(r)\bigr)^{1-\frac{\gamma}{m-1}}, \qquad r\in(s_0,\infty),$$
we can see that $v$ is a nonnegative and nontrivial solution of the following problem
\begin{equation} \label{P_bis}
   \begin{cases}
      -(r^{N+\alpha-1}|v'(r)|^{m-2}v'(r))'= \lambda r^{N+\beta-1}v^{\frac{p(m-1)}{m-1-\gamma}}(r),& r\in (s_0,\infty) \\v(s_0)=1, \qquad v'(s_0)=d\biggl(1-\dfrac{\gamma}{m-1}\biggr) u'(s_0),
   \end{cases}
\end{equation}
where 
\begin{equation}\label{lambda_assunta}
\lambda=\biggl(\frac{m-1-\gamma}{m-1}\biggr)^{m-1}d^{m-1-p}.
\end{equation}
Moreover, $v'(s_0)<0$ by $(H_1)$ and being $u'(s_0)<0$.
Now we apply Lemma \ref{lemma2.2} to problem  \eqref{P_bis} with 
\begin{equation}\label{def_ptilde}
\wp=\frac{p(m-1)}{m-1-\gamma}   
\end{equation}
so that $\wp>m-1$ since $p>m-1-\gamma$ which holds by $p>m-1$ in $(H_1)$.
%Furthermore, \eqref{tesi1_Lemma2.2} in Lemma \ref{lemma2.2} gives
%\begin{equation}\label{tesi1_ptilde}
 %   r^{N+\beta}v^{\frac{(m-1)(p+1)-\gamma}{m-1-\gamma}}(r)\leq Cr^{-\frac{m(m-1-\gamma)(\beta+N)-(N+\alpha-m)[(m-1)(p+1)-\gamma]}{(p+1)(m-1)-\gamma-m(m-1-\gamma)}} 
%\end{equation}
%for all $r\in(s_0,\infty)$, and 
In particular, the identity \eqref{tesi2_Lemma2.2} becomes
\begin{equation}\begin{aligned}\label{equazione}
    &\lambda\left(\dfrac{(N+\beta)(m-1-\gamma)}{(p+1)(m-1)-\gamma}-\dfrac{N+\alpha-m}{m}\right)\int_{s_0}^{r}{s^{N+\beta-1}v^{\frac{(m-1)(p+1)-\gamma}{m-1-\gamma}}(s)ds}\\
    &=-\frac{m-1}m r^{N+\alpha-1}|v'(r)|^{m-1}
   U_{\varrho,v}(r)+\frac{m-1}m s_0^{N+\alpha-1}|v'(s_0)|^{m-1} U_{\varrho,v}(s_0)\\&\quad+
    \dfrac{\lambda(m-1-\gamma)}{(m-1)(p+1)-\gamma}[r^{N+\beta}v^{\frac{(m-1)(p+1)-\gamma}{m-1-\gamma}}(r)-s_0^{N+\beta}],
\end{aligned}\end{equation}
where $U_{\varrho,v}(t)=tv'(t)+(N+\alpha-m)/(m-1)v(t)$.
We study the sign of the left hand side of \eqref{equazione} and we observe that 
$$
   \dfrac{(N+\beta)(m-1-\gamma)}{(p+1)(m-1)-\gamma}-\dfrac{N+\alpha-m}{m}> 0
$$
is valid if  $p$ verifies the following condition
$$
p+1<\frac{m(N+\beta)}{N+\alpha-m}\biggl(1-\frac{\gamma}{m-1}\biggr)+\frac{\gamma}{m-1}    
$$
which is $(H_1)$.
Consequently, since $v$ is positive, we have that  \begin{equation}\label{contradiction_disug}
\lambda\left(\dfrac{(N+\beta)(m-1-\gamma)}{(p+1)(m-1)-\gamma}-\dfrac{N+\alpha-m}{m}\right)\int_{s_0}^{r}{s^{N+\beta-1}v^{\frac{(m-1)(p+1)-\gamma}{m-1-\gamma}}(s)ds}>0\end{equation}
for all $r\in(s_0,\infty)$.
Now, we investigate the sign of the right hand side of 
\eqref{equazione}. First note by Lemma \ref{lemma2.1}, applied to $v$ with $\wp$ given in \eqref{def_ptilde} and $U_\varrho$ replaced with $U_{\varrho,v}$,  we obtain
\begin{equation} \label{segno_primo_termine}
    \frac{m-1}m r^{N+\alpha-1}|v'(r)|^{m-1}
   U_{\varrho,v}(r)\ge0
\end{equation}
for all $r\in(s_0,\infty)$. 
Now we define 
$$L_{s_0}:= \frac{m-1}m s_0^{N+\alpha-1}|v'(s_0)|^{m-1} U_{\varrho,v}(s_0)-\dfrac{\lambda(m-1-\gamma)}{(m-1)(p+1)-\gamma}s_0^{N+\beta},$$ 
so that \eqref{equazione} becomes
\begin{equation}\begin{aligned}\label{costante_L}
    &\lambda\left(\dfrac{(N+\beta)(m-1-\gamma)}{(p+1)(m-1)-\gamma}-\dfrac{N+\alpha-m}{m}\right)\int_{s_0}^{r}{s^{N+\beta-1}v^{\frac{(m-1)(p+1)-\gamma}{m-1-\gamma}}(s)ds}\\
    &\qquad=-\frac{m-1}m r^{N+\alpha-1}|v'(r)|^{m-1}
   U_{\varrho,v}(r)+L_{s_0}\\
    &\qquad \quad+
    \dfrac{\lambda(m-1-\gamma)}{(m-1)(p+1)-\gamma}r^{N+\beta}v^{\frac{(m-1)(p+1)-\gamma}{m-1-\gamma}}(r).
\end{aligned}\end{equation}
In what follows we study the sign of $L_{s_0}$.
By returning to the variable $u$, being $v'(s_0)=d(m-1-\gamma) u'(s_0)/(m-1)$, $u(s_0)=d^{-1}$ and $\lambda$ as in \eqref{lambda_assunta}, we obtain
\begin{equation}\begin{aligned} \label{L}
 L_{s_0}
&=d^{m-1}\left(\frac{m-1-\gamma}{m-1}\right)^{m-1}\biggl\{ \frac{m-1-\gamma}{m}ds_0^{N+\alpha-1}|u'(s_0)|^{m-1}\\
&\hspace{2cm}\cdot\left(s_0u'(s_0)+\frac{N+\alpha-m}{m-1-\gamma} u(s_0)\right)-\dfrac{(m-1-\gamma)d^{-p}s_0^{N+\beta}}{(m-1)(p+1)-\gamma}\biggr\}\\
&=d^{m-1}\left(\frac{m-1-\gamma}{m-1}\right)^{m-1}\biggl\{ \frac{m-1-\gamma}{m}ds_0^{N+\alpha-1}|u'(s_0)|^{m-1}\cdot U_\varrho(s_0)\\
&\hspace{2cm}+\frac{\gamma}{m}\varrho ds_0^{N+\alpha-1}|u'(s_0)|^{m-1}u(s_0)-\dfrac{(m-1-\gamma)d^{-p}s_0^{N+\beta}}{(m-1)(p+1)-\gamma}\biggr\},
\end{aligned} \end{equation}
where in the last equality we have inserted $U_\varrho(s_0)$, defined in Lemma \ref{lemma2.1}.
Furthermore, multiplying the first equation in \eqref{P} by $ru'(r)$ and integrating from $0$ to $s_0$, we have
\begin{equation}\label{int_parti}
-\int_{0}^{s_0}(r^{N+\alpha-1}|u'(r)|^{m-2}u'(r))'ru'(r)dr=\int_{0}^{s_0} r^{N+\beta}u^{p}(r)u'(r)dr.    
\end{equation}
Integrating by parts twice the left hand side of \eqref{int_parti} we obtain
\begin{equation}\label{parti_1}
\begin{aligned}
-\int_{0}^{s_0}&(r^{N+\alpha-1}|u'(r)|^{m-2}u'(r))'ru'(r)dr\\
&=-s_0^{N+\alpha}|u'(s_0)|^m+\int_{0}^{s_0}r^{N+\alpha-1}|u'(r)|^m dr+\frac{1}{m}s_0^{N+\alpha}|u'(s_0)|^m\\
&\hspace{4,5 cm}-\frac{N+\alpha}{m}\int_{0}^{s_0}r^{N+\alpha-1}|u'(r)|^m dr\\
&=-\frac{m-1}{m}s_0^{N+\alpha}|u'(s_0)|^m-\frac{N+\alpha-m}{m}\int_{0}^{s_0}r^{N+\alpha-1}|u'(r)|^m dr.
\end{aligned}
\end{equation}
Now integrating once the right hand side of \eqref{int_parti} we have
\begin{equation} \label{parti_2}
\begin{aligned}
&\int_{0}^{s_0} r^{N+\beta}u^{p}(r)u'(r)dr=\frac{s_0^{N+\beta}u^{p+1}(s_0)}{p+1}-\frac{N+\beta}{p+1}\int_{0}^{s_0}r^{N+\beta-1}u^{p+1}(r)dr.
\end{aligned}
\end{equation}
Using \eqref{parti_1} and \eqref{parti_2} in \eqref{int_parti}, we get
\begin{equation}\label{prima_integrazione}
\begin{aligned}
&-\frac{m-1}ms_0^{N+\alpha}|u'(s_0)|^m-\frac{N+\alpha-m}{m}\int_{0}^{s_0}r^{N+\alpha-1}|u'(r)|^m dr\\
&\quad\quad= \frac{s_0^{N+\beta}u^{p+1}(s_0)}{p+1}-\frac{N+\beta}{p+1}\int_{0}^{s_0}r^{N+\beta-1}u^{p+1}(r)dr.
\end{aligned}
\end{equation}
Moreover, we observe that multiplying  the first equation in \eqref{P} by $u$ and integrating again by parts from $0$ to $s_0$ the left hand side, we have
\begin{equation}\label{seconda_integrazione}
\int_{0}^{s_0}r^{N+\alpha-1}|u'(r)|^mdr=\int_{0}^{s_0}r^{N+\beta-1}u^{p+1}(r)dr-s_0^{N+\alpha-1}|u'(s_0)|^{m-1}u(s_0).
\end{equation}
Inserting \eqref{seconda_integrazione} in \eqref{prima_integrazione}, we see
\begin{equation}\label{eq_finale}
\begin{aligned}
&\frac{m-1}{m}|u'(s_0)|^{m-1}s_0^{N+\alpha-1}U_\varrho(s_0)=\frac{s_0^{N+\beta}u^{p+1}(s_0)}{p+1}\\
&\hspace{3 cm}-\left(\frac{N+\beta}{p+1}-\frac{N+\alpha-m}{m}\right)\int_{0}^{s_0}r^{N+\beta-1}u^{p+1}(r)dr.
\end{aligned}
\end{equation}
Using \eqref{eq_finale} in  \eqref{L}, we obtain
\begin{equation} \label{L_finale}
\begin{aligned}
L_{s_0}&=d^{m-1}\left(\frac{m-1-\gamma}{m-1}\right)^{m-1}\biggl\{ \frac{m-1-\gamma}{m-1}d\biggl[\frac{s_0^{N+\beta}u^{p+1}(s_0)}{p+1}\\
&\quad-\left(\frac{N+\beta}{p+1}-\frac{N+\alpha-m}{m}\right)\int_{0}^{s_0}r^{N+\beta-1}u^{p+1}(r)dr\biggr]\\
&\quad+\frac{\gamma}{m}\varrho ds_0^{N+\alpha-1}|u'(s_0)|^{m-1}u(s_0)-
\dfrac{(m-1-\gamma)d^{-p}s_0^{N+\beta}}{(m-1)(p+1)-\gamma}\biggr\}.
\end{aligned}
\end{equation}
In addition, multiplying by $u(s_0)$ the first equation in \eqref{P} and integrating from $0$ to $s_0$, we have
\begin{equation}\label{eq_maggiorare}\begin{aligned}
|u'(s_0)|^{m-1}u(s_0)s_0^{N+\alpha-1}&=\int_{0}^{s_0}r^{N+\beta-1}u^{p}(r)u(s_0)dr
<\int_{0}^{s_0}r^{N+\beta-1}u^{p+1}(r)dr,
\end{aligned}\end{equation}
where in the last inequality we have used that $u$ is a positive and strictly decreasing solution of \eqref{P}, 
so that $u(s_0)< u(r)$ for all $r\in(0,s_0)$. 

In turn,  
using \eqref{eq_maggiorare} in \eqref{L_finale} and $u(s_0)=d^{-1}$, we obtain
\begin{equation} \label{L_negativo}
\begin{aligned}
\frac{L_{s_0}}{d^{m-1}}&\left(\frac{m-1}{m-1-\gamma}\right)^{m-1}<
\frac{m-1-\gamma}{m-1}\cdot\dfrac{s_0^{N+\beta}u^{p+1}(s_0)}{p+1}d\\
&\quad-\frac{m-1-\gamma}{m-1}\left(\frac{N+\beta}{p+1}-
\frac{N+\alpha-m}{m}\right)d\int_{0}^{s_0}r^{N+\beta-1}u^{p+1}(r)dr\\
&\quad+\frac{\gamma}{m}\varrho d\int_{0}^{s_0}r^{N+\beta-1}u^{p+1}(r)dr-
\dfrac{(m-1-\gamma)d^{-p}s_0^{N+\beta}}{(m-1)(p+1)-\gamma}\\
&=\dfrac{m-1-\gamma}{m-1}s_0^{N+\beta}\left(\frac{du^{p+1}(s_0)}{p+1}-\dfrac{d^{-p}}{p+1-\frac{\gamma}{m-1}}\right)\\
&\quad-d\biggl[\frac{(m-1-\gamma)(N+\beta)}{(m-1)(p+1)}
-\frac{N+\alpha-m}m\biggr]\int_{0}^{s_0}r^{N+\beta-1}u^{p+1}(r)dy\\
&=-\frac{\gamma(m-1-\gamma)}{d^p(m-1)^2(p+1)\left(p+1-\frac{\gamma}{m-1}\right)}s_0^{N+\beta}\\
&\quad-d\left[\frac{(m-1-\gamma)(N+\beta)}{(m-1)(p+1)}-\frac{N+\alpha-m}{m}\right]\int_{0}^{s_0}r^{N+\beta-1}u^{p+1}(r)dy.
\end{aligned}
\end{equation}
Then, since $u$ is positive in $(0,s_0)$, by conditions $(H_0)$, $(H_1)$ and from \eqref{L_negativo}, we have $L_{s_0}<0$. 

Using \eqref{segno_primo_termine} and the negativity of $L_{s_0}$ in \eqref{costante_L} we have
\begin{equation} \label{ineq_limite}
\begin{aligned}
&\lambda\left[\dfrac{(N+\beta)(m-1-\gamma)}{(p+1)(m-1)-\gamma}-\dfrac{N+\alpha-m}{m}\right]\int_{s_0}^{r}{s^{N+\beta-1}v^{\frac{(m-1)(p+1)-\gamma}{m-1-\gamma}}(s)ds}\\
&\hspace{2 cm}\le\dfrac{\lambda(m-1-\gamma)}{(m-1)(p+1)-\gamma}r^{N+\beta}v^{\frac{(m-1)(p+1)-\gamma}{m-1-\gamma}}(r).
\end{aligned}
\end{equation}
Finally, we observe that, since $\wp$ given in \eqref{def_ptilde} satisfies $(H_1)'$, which is in force since $p$ verifies $(H_1)$, then Corollary \ref{cor_p} holds and we obtain \begin{equation} \label{infinitesimo_tesi1}
    r^{N+\beta}v(r)^{\frac{(m-1)(p+1)-\gamma}{m-1-\gamma}}\to0 \qquad \text{as} \; r\to\infty.
\end{equation}
Thus, by letting $r\to\infty$ in \eqref{ineq_limite}, since \eqref{infinitesimo_tesi1} holds, we get $$\lambda\left[\dfrac{(N+\beta)(m-1-\gamma)}{(p+1)(m-1)-\gamma}-\dfrac{N+\alpha-m}{m}\right]\int_{s_0}^{\infty}{s^{N+\beta-1}v^{\frac{(m-1)(p+1)-\gamma}{m-1-\gamma}}(s)ds}\le0,$$
so that, by letting $r\to\infty$ in \eqref{contradiction_disug}, we arrive to
$$\int_{s_0}^{\infty}{s^{N+\beta-1}v^{\frac{(m-1)(p+1)-\gamma}{m-1-\gamma}}(s)ds}=0.$$
%yielding $v=u\equiv0$ in $(s_0,\infty)$ thanks to \eqref{cambio_var}. 
From this contradiction,  the conclusion of the proof follows.
\end{proof}

\section{A priori estimates}\label{stima}
In this section we use the Liouville type results contained in the previous section to get a priori estimates for positive radial solutions of problem \eqref{PROB}. Thus, we introduce, for each $k\in\mathbb{N}$,  the following function
$$
 g_k(s):=(g\circ T_k)(s)=g(T_k(s)),
 \quad s\ge 0,
$$
where  $T_k(s):=\max\{-k,\min\{k,s\}\}$,  $s\in\mathbb R$ is the well known truncated function.
In particular, we have 
$$g_k(s)=g(k) \text{ if } s\ge k,\quad  g_k(s)=g(s) 
\text{ if } 0\le s<k.$$ 
Now we consider the family of truncated problems parametrized by $\xi\ge0$,
\begin{equation} \label{prob_parametrizzato}
{\small\begin{cases}
  -\!\left(\!\dfrac{r^{N+\alpha-1}|v'(r)|^{m-2}v'(r)}{\big(a(r)+g_k(v(r))\big)^\gamma}\!\right)'=r^{N+\beta-1}\left(v^p(r)+\dfrac{\xi}{h(\|v\|_\infty)}\right), \quad\,\;0<r<R,
  \\v'(0)=0, \qquad v(R)=0,
   \end{cases}}
\end{equation}
where $h$ is defined as follows 
\begin{equation} \label{def_h}
h(t)=\begin{cases} t^{q-p} \quad &\text{if}\quad t\ge1,\\1 \quad &\text{if}\quad 0\le t\le1, \end{cases}    
\end{equation}
with $q>p$.
%\begin{figure}
%\centering
%\fbox{\includegraphics[width=.5\textwidth, height=.3\textheight, keepaspectratio]{Grafico troncata.png}}
%\caption{Truncated function ($T_k$)}\label{troncata}
%\end{figure}

We give now a crucial result, based on the celebrated blow up technique in \cite{GS}, which provides a priori bound for solutions of the truncated problem \eqref{prob_parametrizzato}.
\begin{theorem}\label{thm3.1}
Assume $(H_0)$, $(H_1)'$, $(H_2)$ and $(H_3)$. Then there is a positive constant $C_k$, which depends only on $k$, such that 
\begin{equation}\label{stimaapriori}
    \|v\|_\infty\le C_k,
\end{equation} for every positive solution $v$ of problem \eqref{prob_parametrizzato}.
\end{theorem}
\begin{proof}
Let be $k\in\mathbb{N}$ be fixed and assume by contradiction that there is a sequence of positive solutions $(v_n)_n$ of problem \eqref{prob_parametrizzato}, such that $\|v_n\|_\infty\to\infty$ for $n\to\infty$.

We observe that, since $v_n>0$ for all $n\in\mathbb N$, $\xi \ge0$ and $h$ is a positive function,  by a qualitative analysis of the equation in \eqref{prob_parametrizzato}, we obtain that each $v_n$ is a decreasing function.
Now we use the following changes of variables 
\begin{equation}\label{3.2}
y=\frac{z_n}{t_n}r,    
\end{equation}
where
\begin{equation}\label{zktk}
t_n:=\|v_n\|_\infty, \qquad z_n:=\big(a(0)+g(k))^{\frac{\gamma}{\beta-\alpha+m}}t_n^{\frac{\beta-\alpha+1+p}{\beta-\alpha+m}}>0,\end{equation}
so that $t_n\to\infty$ as $n\to\infty$. Define 
\begin{equation}\label{defwn}
w_n(y):=\frac{v_n(r)}{t_n},
\end{equation}
clearly $\|w_n\|_\infty=1$. Since $v_n$ is a positive and decreasing solution of \eqref{prob_parametrizzato} for all $n\in\mathbb{N}$, we observe that $w_n(0)=1$ by virtue of the definition of $\|\cdot\|_\infty$, then for  large $n$, using the definition of $h$ given in \eqref{def_h} and the fact that $\|w_n\|_\infty=1$, the function $w_n$ is a solution of the following problem
\begin{equation} \label{prob_w_n}
  {\footnotesize\begin{cases}
 -\!\left(\!\dfrac{y^{N+\alpha-1}|w_n'(y)|^{m-2}w_n'(y)}{\big(a(t_nz_n^{-1}y)+g_k(t_nw_n(y))\big)^\gamma}\!\right)'\!=\!\dfrac{t_n^{\beta-\alpha+1}}{z_n^{\beta-\alpha+m}}y^{N+\beta-1}\!\biggl(t_n^pw_n^p(y)+\dfrac{\xi}{t_n^{q-p}}\biggr)\!, \quad y<\dfrac{z_nR}{t_n},
  \\w_n'(0)=0, \qquad w_n(0)=1, \qquad w_n\biggl(\dfrac{Rz_n}{t_n}\biggr)=0,
   \end{cases}}
\end{equation}
where $'$ denotes the derivative respect to the variable of the function under consideration, hence  
$w_n'(y)=\frac 1{z_n}v_n'(r)$
yielding $w_n'(0)=0$, $w_n'<0$ in $(0,z_nR/t_n)$ and
$$\biggl (|v_n'(r)|^{m-2}v_n'(r) \biggr)'
%=\frac d{dr}\biggl (|v_n'(r)|^{m-2}v_n'(r) \biggr)\\&
%=\frac{z_n^m}{t_n}\frac d{dy}\biggl (|w_n'(y)|^{m-2}w_n'(y) \biggr)
=\frac{z_n^m}{t_n}\biggl (|w_n'(y)|^{m-2}w_n'(y) \biggr)',$$
so that
$$\biggl (\frac{r^{N+\alpha -1}|v_n'(r)|^{m-2}v_n'(r)}{\big(a(r)+g_k(v_n(r))\big)^\gamma} \biggr)'
=\frac{t_n^{N+\alpha-2}}{z_n^{N+\alpha -1-m}}
\biggl (\frac{y^{N+\alpha -1}|w_n'(y)|^{m-2}w_n'(y)}{\big(a(t_nz_n^{-1}y)+g_k(t_nw_n(y))\big)^\gamma} \biggr)',$$
in turn equation in \eqref{prob_w_n} follows immediately from \eqref{prob_parametrizzato}.
We also note that 
\begin{equation} \label{R+}
\lim_{n\to\infty}\frac{Rz_n}{t_n}=\lim_{n\to\infty} R\big(a(0)+g(k)\big)^{\frac{\gamma}{\beta-\alpha+m}}t_n^{\frac{1+p-m}{\beta-\alpha+m}}=\infty,
\end{equation}
since $(1+p-m)/(\beta-\alpha+m)>0$ by $(H_0)$ and $(H_1)$, and $t_n\to\infty$ as $n\to\infty$ by contradiction.

Observe that integrating the equation in \eqref{prob_w_n} from $0$ to $y\in\left(0,z_nR/t_n\right)$ and replacing $z_n$ in the right hand side, we have 
\begin{equation}\label{integrazione_wn}
\frac{y^{N+\alpha-1}|w_n'(y)|^{m-1}}{\big(a(t_nz_n^{-1}y)+g_k(t_nw_n(y))\big)^\gamma}=\frac{\int_{0}^{y}\tau^{N+\beta-1}\left(w_n^p(\tau)+\frac{\xi}{t_n^{q}}\right)d\tau}{\bigl(a(0)+g(k)\bigr)^\gamma}.
\end{equation}
Since $w_n$ is a positive and decreasing solution of \eqref{prob_w_n} for large $n$ and being $t_n\to\infty$ as $n\to\infty$ then, for large $n$, we have $w_n^p(\tau)+\xi/t_n^q\le w_n^p(0)+1=2$ for all parameter $\xi\in\mathbb R^+$, thus \eqref{integrazione_wn} gives
\begin{equation}\label{lip_1}
|w_n'(y)|^{m-1}\le \frac{2}{N+\beta}\biggl(\frac{a(t_nz_n^{-1}y)+g_k(t_nw_n(y))}{a(0)+g(k)}\biggr)^\gamma y^{\beta-\alpha+1}.   
\end{equation}
Now using $(H_3)$ and that $g_k(s)\le g(k)$ for all $s\ge0$ we obtain for all $n\in\mathbb N$
\begin{equation}\label{frac<1}
0\le \biggl(\frac{a(t_nz_n^{-1}y)+g_k(t_nw_n(y))}{a(0)+g(k)}\biggr)^\gamma\le \biggl(\frac{c_2+g(k)}{a(0)+g(k)}\biggr)^\gamma\le \biggl(1+\frac{c_2}{c_1}\biggr)^\gamma=:C_\gamma,\end{equation}
thus,  \eqref{lip_1} combined with \eqref{frac<1} gives
for all $n$,
 \begin{equation} \label{limitatezza_derivata}
|w_n'(y)|\le \biggl(\frac{2 C_\gamma}{N+\beta}y^{\beta-\alpha+1}\biggr)^{\frac{1}{m-1}}.
 \end{equation}
From \eqref{limitatezza_derivata} and $\beta-\alpha+1>0$ by $(H_0)$, we get that $w_n'$ is uniformly bounded in compact intervals.

Let now $\bar R$ be a positive number such that $\bar R<Rz_n/t_n$ for large $n$ and we consider the restriction of $w_n$ to $[0,\bar R]$, stille denoted with $w_n$. Then, there is a constant $C(\bar R)>0$ so that $
|w_n'(y)|\le C(\bar R)$ for all $n\in\mathbb{N}$ and for all $y\in[0,\bar R]$, in turn
the sequence $(w_n)_n$ is uniformly equilipschitz or, equivalently, uniformly equicontinuous. 

Furthermore, from $\|w_n\|_\infty=1$, the sequence $(w_n)_n$ is also uniformly bounded in $[0,\bar R]$. By Ascoli Arz\`ela's Theorem, $(w_n)_n$ contains a subsequence converging uniformly  (which we still denote by $(w_n)_n$), namely 
\begin{equation} \label{converg_wn}
 w_n\to w \quad\text{in}\;\;C[0,\bar R].   
\end{equation}
We claim that, for $\bar R$ fixed, if we define $w_n(\bar R)=\delta_n$, then necessarily there exists $\bar \varepsilon$ such that 
\begin{equation} \label{w_n(R)limit}
0<\bar \varepsilon\le\delta_n<1\quad \text{for all}\quad n\in\mathbb N.    
\end{equation}
Indeed, since $w_n$ is a decreasing function for all $n\in\mathbb N$ we have that $0<\delta_n<1$ and if we assume by contradiction that $\displaystyle\inf_n\delta_n=0$ then for $\varepsilon>0$ fixed there exists $\bar n$ sufficiently large such that 
\begin{equation}\label{contradict_equilip}w_{\bar n}(\bar R)=\delta_{\bar n}<\frac\varepsilon 2.\end{equation}
On the other hand, as we have observed, $w_n$ is a uniformly equicontinuous function in $[0,\bar R]$, then for all $\varepsilon>0$ fixed there exists $\eta(\varepsilon)=\eta$ such that for all $y_1$ and $y_2$ with 
\begin{equation} \label{defunifequicont}
\quad |y_1-y_2|<\eta \quad\text{then} \quad |w_{n}(y_1)-w_{n}(y_2)|<\frac{\varepsilon}{2}   
\end{equation}
for all $n\in\mathbb N$. 
Let $\bar y<\bar R$ such that $0<\bar R-\bar y<\eta$ then applying \eqref{defunifequicont} with $y_1=\bar y$,  $y_2=\bar R$ and $n=\bar n$,   since $w_{\bar n}'<0$,  we have
\begin{equation}\label{epsilon2}
  0<w_{\bar n}(\bar y)-w_{\bar n}(\bar R)<\frac{\varepsilon}{2}.
\end{equation}
Now adding \eqref{contradict_equilip} and \eqref{epsilon2}
we arrive to 
\begin{equation}\label{cointr_2}w_{\bar n}(\bar y)<\varepsilon.\end{equation}
Take $\tilde y$ such that $0<\tilde y-\bar y<\eta$,
so that, by \eqref{defunifequicont} with $y_1=\tilde y$ and $y_2=\bar y$, it holds
$$0<w_{\bar n}(\tilde y)-w_{\bar n}(\bar y)<\frac{\varepsilon}{2},$$
which, added to \eqref{cointr_2}, gives
$w_{\bar n}(\tilde y)<\frac32\varepsilon.$
Iterating this procedure, we find a point $y^*$ sufficiently close to zero, such that, being  $w_{\bar n}(0)=1$, we have
$$\frac 12 
\le w_{\bar n}(y^*)<C\varepsilon, \quad C>0.$$
The arbitrariness of $\varepsilon$, concludes the proof of \eqref{w_n(R)limit}.

Hence  \eqref{w_n(R)limit} gives $w_n(\bar R)\ge C>0$ for all $n$, from which it we immediately follows 
\begin{equation} \label{wpositiva}
    w(y)>0 \quad\text{for all} \quad y\in[0,\bar R]
\end{equation}
Since \eqref{wpositiva} holds and by $t_n\to\infty$ for $n\to\infty$ we see that $t_nw_n(s)\to\infty$ for $n\to\infty$ for all $s\in[0,\bar R]$ fixed, then for $n$ sufficiently large we have that \begin{equation} \label{opiccoloxi}
g_k(t_nw_n(y))=g(k) \quad \text{for all} \quad y\in[0,\bar R]    
\end{equation} 
and 
\begin{equation}\label{refff}
\xi/t_n^q=o(1)\quad \text{as} \quad n\to\infty,
\end{equation}
for any $\xi\in\mathbb R^+$.
Moreover we observe that
\begin{equation} \label{argom_a}
\frac{t_n}{z_n}=\frac{t_n^{(m-1-p)/(\beta-\alpha+m)}}{\big(a(0)+g(k)\big)^{\gamma/(\beta-\alpha+m)}}\to0\quad\text{as}\,\,n\to\infty
\end{equation}
by $(H_0)$ and $(H_1)$ and $t_n\to\infty$ as $n\to\infty$ by contradiction.
Then, by continuity of $a$ assumed in $(H_3)$ and by \eqref{argom_a}, we have
\begin{equation} \label{a(0)}
 a(t_nz_n^{-1}y)\to a(0) \quad\text{as} \quad n\to\infty,
\end{equation}
for all $y\in[0,\bar R]$. In particular, from \eqref{a(0)} and for $y\in[0,\bar R]$, we obtain
\begin{equation} \label{limite_frac}
\lim_{n\to\infty}\biggl(\frac{a(t_nz_n^{-1}y)+g(k)}{a(0)+g(k)}\biggr)^\gamma=1.
\end{equation}
Using \eqref{opiccoloxi}, \eqref{limite_frac}, \eqref{refff} and \eqref{integrazione_wn} we get that for $y\in[0,\bar R]$ and $n$ sufficiently large is valid
$$|w_n'(y)|^{m-1}=\frac{1+o(1)}{y^{N+\alpha-1}} \int_{0}^{y}\tau^{N+\beta-1}\big(w_n^p(\tau)+o(1)\big)d\tau,$$
then $\bar R<z_nR/t_n$.
Consequently by integration afrom $0$ to $y\in[0,\bar R]$ to obtain, for $n$ sufficiently large, we arrive to
\begin{equation}\label{w_n_equazione}
1-w_n(y)=[1+o(1)]^{\frac{1}{m-1}}\int_{0}^{y}\biggl[\frac{1}{s^{N+\alpha-1}}\int_{0}^{s}\tau^{N+\beta-1}\big(w_n^p(\tau)+o(1)\big)d\tau\biggr]^{\frac{1}{m-1}}ds    
\end{equation}
for all $y\in[0,\bar R]$.
Passing to the limit for $n\to\infty$ in \eqref{w_n_equazione}, by using \eqref{converg_wn} and Lebesgue's dominated Theorem we obtain that $w$ satisfies the following integral equation
\begin{equation} \label{eq_integrale_w}
1-w(y)=\int_{0}^{y}\biggl[\frac{1}{s^{N+\alpha-1}}\int_{0}^{s}\tau^{N+\beta-1}w^p(\tau)d\tau\biggr]^{\frac{1}{m-1}}ds.
\end{equation}
Furthermore, \eqref{eq_integrale_w} gives by differentiation, that $w$ is a decreasing and $C^1[0,\bar R]$ function with  $w'(0)=0$ indeed, by continuity of $w'$ we get
\begin{equation} \label{w'=0}
|w'(0)|^{m-1}\le\lim_{y\to0^+}\frac{w^p(0)}{y^{N+\alpha-1}}\int_{0}^{y}\tau^{N+\beta-1}d\tau=\lim_{y\to0^+}\frac{y^{\beta-\alpha+1}}{N+\beta}=0,
\end{equation}
where we have used that $w$ is decreasing and  $w(0)=1$.

In particular, being $(w_n)_n$ a sequence of positive and decreasing functions then its uniform limit $w$, satisfying \eqref{eq_integrale_w},  is a positive solution in $[0,\bar R]$ of the following problem
\begin{equation} \label{problema_w}
\begin{cases}
  -\big(y^{N+\alpha-1}|w'(y)|^{m-2}w'(y)\big)'=y^{N+\beta-1}w^p(y), \quad y\in(0,\bar R),
  \\w(0)=1, \qquad w'(0)=0.
\end{cases}
\end{equation}
Arguing as in the proof of Proposition 4.1 in \cite{CMM}, we claim that $w$ can be extended to the entire $\mathbb{R}^+$, obtaining a positive solution of \eqref{p2} with $\lambda=1$.
To see this, it is sufficient to note that by \eqref{R+}  we can repeat the above argument on an interval $[0,R^*]$, with $R^*>\bar R$ for the convergent sequence (on $[0,\bar R]$) $(w_n)_n$. In this manner we obtain a function $\tilde w$, solution of  \eqref{problema_w} on $[0,R^*]$, that satisfies
$$w(y)=\tilde w(y)$$
in $[0,\bar R]$. It is now clear that $w$ can be extended to $\mathbb R^+$ as a positive solution of problem \eqref{problema_w} in  $(0,\infty)$ and the claim follows. 
Thus, $w$ is a positive solution of \eqref{p2}, contradicting  Theorem \ref{teorema1} applied with $\wp=p$ and $\lambda=1$.

In turn, the a priori estimate \eqref{stimaapriori} holds with $C_k$ independent on $\xi$.
\end{proof}

\section{An existence result for the truncated problem}\label{probtronc}

In this section we will give an existence result related to the truncated problem \eqref{prob_parametrizzato} with $\xi=0$, namely
\begin{equation} \label{xi0}
\begin{cases}
  -\left(\dfrac{r^{N+\alpha-1}|v'(r)|^{m-2}v'(r)}{\big(a(r)+g_k(v(r))\big)^\gamma}\right)'=r^{N+\beta-1}v^p(r), \quad 0<r<R,
  \\v'(0)=0, \qquad v(R)=0.
   \end{cases}
\end{equation}

The proof of the existence of positive solutions of problem \eqref{xi0} is based on Theorem \ref{kran_th} and in order to use it, we consider the Banach space $X=C[0,R]$ endowed with the $L^\infty$-norm, and the following cone  by $\mathcal{C}=\{v\in C[0,R];\quad v\ge0,\;v(R)=0\}$.

Define the operator $F:X\to X$ by
\begin{equation}\label{operatore_F}
F(v)(r)=\int_{r}^{R}\biggl[\frac{\bigl(a(s)+g_k(v(s))\bigr)^\gamma}{s^{N+\alpha-1}}\int_{0}^{s}\tau^{N+\beta-1}|v(\tau)|^p d\tau\biggr]^{\frac{1}{m-1}}ds.
\end{equation}
Note that every fixed points $v$ of the operator $F$  are positive solutions of problem \eqref{prob_parametrizzato} with $\xi=0$, namely $v$ is a solution of problem \eqref{xi0}.
\begin{lemma}
The operator $F:X\to X$ defined in \eqref{operatore_F} is compact, and the cone $\mathcal{C}$ is invariant under $F$, that is $F(\mathcal{C})\subset\mathcal{C}$.
\end{lemma}
\begin{proof}
In order to prove that $F$ is a compact operator, we claim that, given a sequence 
\begin{equation}\label{vn_bounded}(v_n)_n\subset X \quad\text{such that}\quad \|v_n\|_\infty\le \overline{C},\end{equation} 
for some $\overline{C}$,  then $(F(v_n))_n$ is equicontinuous and uniformly bounded in $X$. Consequently, using the Ascoli Arz\`ela's Theorem, we have that $(F(v_n))_n$ converges in $C[0,R]$, so that the compactness of $F$ is proved. 
To reach the claim, first we observe that the sequence $(F(v_n))_n$ is uniformly bounded, indeed
$$
|F(v_n)(r)|\le\int_{r}^{R}\frac{\bigl(a(s)+g_k(v(s))\bigr)^\frac{\gamma}{m-1}}{s^{\frac{N+\alpha-1}{m-1}}}
\biggl(\int_{0}^{s}\tau^{N+\beta-1}|v_n(\tau)|^p
d\tau\biggr)^{\frac{1}{m-1}}ds\\
$$
then, using the boundedness of $a$ given in $(H_3)$ and $g_k(v)\le g(k)$ for every $v\ge 0$, by the definition of $g_k$, we have
\begin{equation} \label{Fv_n}\begin{aligned}
|F(v_n)(r)|&\le\bigl(c_2+g(k)\bigr)^{\frac{\gamma}{(m-1)}}\|v_n\|_\infty^{\frac{p}{m-1}}
\int_{r}^{R}\frac1{s^{\frac{N+\alpha-1}{m-1}}} \biggl(\int_{0}^{s}\tau^{N+\beta-1}d\tau\biggr)^{\frac{1}{m-1}}ds
\\&\le\biggl[\frac{\bigl(c_2+g(k)\bigr)^\gamma}{N+\beta}\biggr]^{\frac{1}{m-1}}\|v_n\|_\infty^{\frac{p}{m-1}}\cdot\frac{m-1}{\beta-\alpha+m}\cdot R^{\frac{\beta-\alpha+m}{m-1}}.    
\end{aligned}\end{equation}
Then using the uniform boundedness of $(v_n)_n$ given in \eqref{vn_bounded}, we obtain 
$$|F(v_n)(r)|\le\overline{C}\biggl[\frac{\bigl(c_2+g(k)\bigr)^\gamma}{N+\beta}\biggr]^{\frac{1}{m-1}}\cdot\frac{m-1}{\beta-\alpha+m}R^{\frac{\beta-\alpha+m}{m-1}}$$
for all $n\in\mathbb N$ and $r\in [0,R]$, namely the sequence $(F(v_n))_n$
is uniformly bounded.

Now we prove that the sequence $(F(v_n))_n$ is equilipschitz, or equivalently equicontinuous, indeed,
 using the same tools as above, we get by $(H_0)$
%$(H_2)$, $g_k(s)\le g(k)$ and $v_n\le\|v_n\|_\infty$
%\begin{align*}
%|(F'v_n)(r)|&=\biggl[\frac{\bigl(a(r)+g_k(v(r))\bigr)^\gamma}{r^{N+\alpha-1}}\int_{0}^{r}\tau^{N+\beta-1}v_n^p(\tau)d\tau\biggr]^{\frac{1}{m-1}}\\
%&\le\biggl[\frac{\bigl(c_2+g(k)\bigr)^\gamma}{N+\beta}\biggr]^{\frac{1}{m-1}}\cdot \|v_n\|^{\frac{p}{m-1}}_\infty\cdot r^{\frac{\beta-\alpha+1}{m-1}}, 
%\end{align*}
%in particular using $(H_0)$
$$
\begin{aligned}
|F'(v_n)(r)|&\le\overline{C}\biggl[\frac{\bigl(c_2+g(k)\bigr)^\gamma}{N+\beta}\biggr]^{\frac{1}{m-1}}\cdot r^{\frac{\beta-\alpha+1}{m-1}}\le C R^{\frac{\beta-\alpha+1}{m-1}},
\end{aligned}
$$
for all $n\in\mathbb{N}$ and $r\in[0,R]$, which leads to our conclusion. Finally,
the invariance of the cone $\mathcal{C}$ under $F$ is due to the positivity of functions in $\mathcal{C}$, the fact that $F(v)(R)=0$ for all $v\in \mathcal{C}$, by the definition of $F$, and the regularity of $F$.
\end{proof}

Now we are ready to prove the existence result of the truncated problem \eqref{prob_parametrizzato} with $\xi=0$, that is problem \eqref{xi0}.
\begin{theorem} \label{th_4.3}
Assume $(H_0)$, $(H_1)'$, $(H_2)$ and $(H_3)$. Then there exists a positive solution of the truncated problem \eqref{xi0}.
\end{theorem}
\begin{proof}
To prove the existence of at least a positive solution for the truncated problem \eqref{xi0}, it is enough to show that $F$ has a fixed point in $\mathcal{C}$. For this claim, we will verify conditions $(a)-(d)$ of Theorem \ref{kran_th} given in Section \ref{cap1}. Define the homotopy $\mathcal H:[0,1]\times\mathcal{C}\to\mathcal{C}$ by
$$
\mathcal H(t,v)(r)=\int_{r}^{R}\biggl[\frac{\bigl(a(s)+g_k(v(s))\bigr)^\gamma}{s^{N+\alpha-1}}\int_{0}^{s}\tau^{N+\beta-1}\biggl(v^p(\tau)+\frac{t\xi}{h(\|v\|_\infty)}\biggr)d\tau\biggr]^{\frac{1}{m-1}}ds,   
$$
with $\xi\ge0$ to be chosen. Note that, similarly for the operator $F$  we can show that $\mathcal H$ is a compact homotopy, being $h\ge1$ and $0\le t\le 1$.
Furthermore,  $\mathcal H(0,v)(r)=F(v)(r)$,
so that $(b)$ is verified.

On the other hand, in order to verify $(a)$
we use that, by \eqref{Fv_n},
\begin{equation} \label{tF(v)}
\begin{aligned}
\|tF(v)\|_\infty&=\sup_{r\in[0,R]}|tF(v)(r)|\le \biggl[\frac{\bigl(c_2+g(k)\bigr)^\gamma}{N+\beta}\biggr]^{\frac{1}{m-1}}\|v\|_\infty^{\frac{p}{m-1}}\cdot\frac{m-1}{\beta-\alpha+m}\cdot R^{\frac{\beta-\alpha+m}{m-1}}
\end{aligned}
\end{equation}
for all $t\in[0,1]$ and $v\in\mathcal C$.
Now, let $v\in \mathcal C\setminus \{0\}$ with $\delta:=\|v\|_\infty\in (0,1)$ sufficiently small  and such that 
$$\biggl[\frac{\bigl(c_2+g(k)\bigr)^\gamma}{N+\beta}\biggr]^{\frac{1}{m-1}}\delta^{\frac{p-m+1}{m-1}}\cdot\frac{m-1}{\beta-\alpha+m}\cdot R^{\frac{\beta-\alpha+m}{m-1}}<1,$$
this is possible since by the definition of $\mathcal C$.
Hence, for all $v$ as above, from \eqref{tF(v)} we have $\|tF(v)\|_\infty<\|v\|_\infty$ for all $t\in[0,1]$, which gives $tF(v)\ne v$ for all $t\in[0,1]$. So that $(a)$ is valid.

To obtain (c), namely there exists $\eta>\delta$ such that $\mathcal H(t,v)\ne v$ for all $\|v\|_\infty=\eta$ and $t\in[0,1]$, it is enough to  choose $\eta>\max\{C_k,1\}$, where $C_k$ is defined in Theorem \ref{thm3.1}. Indeed, if we take $\eta>\delta$ being $\delta$ small by $(a)$, then necessarily $(c)$ holds since
if $\mathcal H(t_\eta,v_\eta)=v_\eta$ for some $v_\eta\in\mathcal{C}$ with $\|v_\eta \|_\infty=\eta$ and $t_\eta\in[0,1]$,  then $v_\eta$ is a positive solution of \eqref{prob_parametrizzato} with $\xi$ replaced by $t_\eta\xi$ and Theorem \ref{thm3.1} gives 
$\eta=\|v_\eta\|_\infty\le C_k<\eta$, which yields a contradiction.

Finally, we check the last condition $(d)$ which requires that for $\eta>\delta$ given in $(c)$, then $\mathcal H(1,v)\ne v$ for all $v\in\mathcal{C}$ with $\|v\|_\infty\le\eta$.
We claim that choosing $\xi$ sufficiently large in the definition of $\mathcal H$, if $\mathcal H(1,v)=v$ for some $v\in\mathcal{C}$, that is $v$ is a solution of \eqref{prob_parametrizzato}, then necessarily
\begin{equation}\label{bound_below}
\|v\|_\infty>\max\{1, C_k\}\ge C_k.
\end{equation}
In particular, let $v\in\mathcal{C}$ such that $\mathcal H(1,v)=v$,  using that $v$ is a nonincreasing positive solution of \eqref{prob_parametrizzato} we have that
\begin{equation}\label{norm_inf}
\begin{aligned}
\|v\|_\infty&=v(0)=\mathcal H(1,v)(0)\\
&=\!\int_{0}^{R}\!\biggl[\frac{\bigl(a(s)+g_k(v(s))\bigr)^\gamma}{s^{N+\alpha-1}}
\!\int_{0}^{s}\tau^{N+\beta-1}\!\biggl(v^p(\tau)+
\frac{\xi}{h(\|v\|_\infty)}\biggr)d\tau\biggr]^{\frac{1}{m-1}}\!ds.
\end{aligned}\end{equation}
From \eqref{norm_inf}, $(H_2)$, $(H_3)$ and $v\ge0 $ we get
\begin{equation} \label{stima_norma}
\begin{aligned}
\|v\|_\infty&\ge\big(c_1+g(0)\big)^{\frac{\gamma}{m-1}}\int_{0}^{R}
\biggl[\frac{1}{s^{N+\alpha-1}}
\int_{0}^{s}\tau^{N+\beta-1}\cdot \frac{\xi}{h(\|v\|_\infty)} d\tau\biggr]^{\frac{1}{m-1}}ds\\
&=\big(c_1+g(0)\big)^{\frac{\gamma}{m-1}}\biggl(\frac{\xi}{(N+\beta)h(\|v\|_\infty)}\biggr)^{\frac{1}{m-1}}\int_{0}^{R}s^{\frac{\beta-\alpha+1}{m-1}}ds=C\biggl(\frac{\xi}{h(\|v\|_\infty)}\biggr)^{\frac{1}{m-1}},
\end{aligned}
\end{equation}
where 
$$C=\frac{(c_1+g(0))^{\frac{\gamma}{m-1}}(m-1)}{\beta-\alpha+m}\cdot
(N+\beta)^{-\frac{1}{m-1}}R^{\frac{\beta-\alpha+m}{m-1}}.$$
Now, because of \eqref{def_h} and by the choice of $\eta$, the case $C_k\le\|v\|_\infty\le\eta$ is not possible by Theorem \ref{thm3.1}. Hence, we have two cases: either $\|v\|_\infty\le\min\{1,C_k\}$ or $1<\|v\|_\infty\le C_k$. In the first case, by \eqref{def_h} and \eqref{stima_norma},
we get 
$$
\|v\|_\infty\ge C \xi^{\frac{1}{m-1}},
$$
thus  choosing $\xi$ in the homotopy $\mathcal H$ sufficiently large, say  $\xi>C^{1-m}$, we obtain
\begin{equation} \label{ge_1}
\|v\|_\infty>1.   
\end{equation} 
Differently, if $1<\|v\|_\infty\le C_k$, using \eqref{def_h}, \eqref{stima_norma} and $p<q$ we obtain $$\|v\|_\infty\ge C \xi^{\frac{1}{m-1}}\|v\|_\infty^{\frac{p-q}{m-1}}\ge C \xi^{\frac{1}{m-1}}C_k^{\frac{p-q}{m-1}},$$
then choosing $\xi$ sufficiently large, say $\xi> C^{1-m}C_k^{m-1-p+q}$,  we get 
\begin{equation} \label{ge_2}\|v\|_\infty>C_k.
\end{equation}
By inequalities \eqref{ge_1} and \eqref{ge_2}, choosing the parameter $\xi$ in the homotopy $\mathcal H$ such that $\xi>\max\{C^{1-m},C^{1-m}C_k^{m-1-p+q}\}$ we obtain \eqref{bound_below}, which gives the contradiction required. Consequently, if we choose $\|v\|_\infty\le\eta$, we have $\mathcal H(1,v)\ne v$ so that the proof of $(d)$ is concluded.

In conclusion, using Theorem \ref{kran_th}, we have that the operator $F$, defined in \eqref{operatore_F}, has a fixed point $v\in \mathcal C$, which is a positive solution of \eqref{xi0} such that 
\begin{equation}\label{norma_limitata}
\delta<\|v\|_\infty\le C_k,
\end{equation}
being $\eta>C_k$.
\end{proof}
\begin{remark}
As a consequence of Theorem \ref{th_4.3} we have a result of the type of \cite[Proposition 3.2]{R}. In particular, it holds the following: there exists $\xi^*>0$ such that problem \eqref{prob_parametrizzato} has no positive solutions for any $\xi\ge\xi^*$. Indeed, let $v$ be a positive solution of \eqref{prob_parametrizzato}, then $v$ satisfies the following integral equation
$$
v(r)=\int_{r}^{R}\!\biggl[\frac{\bigl(a(s)+g_k(v(s))\bigr)^\gamma}{s^{N+\alpha-1}}
\!\int_{0}^{s}\tau^{N+\beta-1}\!\biggl(v^p(\tau)+
\frac{\xi}{h(\|v\|_\infty)}\biggr)\!d\tau\biggr]^{\frac{1}{m-1}}\!ds.
$$
Thus, by $v'(r)<0$ we get
$$
\|v\|_\infty=\int_{0}^{R}\!\biggl[\frac{\bigl(a(s)+g_k(v(s))\bigr)^\gamma}{s^{N+\alpha-1}}
\!\int_{0}^{s}\tau^{N+\beta-1}\!\biggl(v^p(\tau)+
\frac{\xi}{h(\|v\|_\infty)}\biggr)\!d\tau\biggr]^{\frac{1}{m-1}}\!ds.
$$
Using the same argument in the proof of the validity of condition $(d)$ in Theorem \ref{th_4.3}, we get that for $\xi\ge\xi^*$, with $\xi^*$ sufficiently large, necessarily $\|v\|_\infty>C_k$, which contradicts the a priori estimate for positive solutions of problem \eqref{prob_parametrizzato} given by Theorem \ref{thm3.1}.
\end{remark}

\section{Proof of Theorem 1}\label{exixtence}

We are now ready to prove the main existence theorem of the paper Theorem \ref{th_1.1}. Note that Theorem \ref{th_4.3} and Theorem \ref{th_1.1} are strictly connected. Indeed positive solutions of \eqref{xi0} with particular values of $k$ are positive radial solutions of \eqref{PROB}.

\begin{proof}[Proof of Theorem \ref{th_1.1}] We claim that there exists  $k_0\in \mathbb N$ such that the 
corresponding  positive solution of the  truncated problem \eqref{xi0} with $k=k_0$, given by Theorem \ref{th_4.3} and
denoted with $v_{k_0}$, verifies
\begin{equation}\label{estimate_k0}\|v_{k_0}\|_\infty\le k_0.
    \end{equation}
Indeed, we first observe that  $\delta<\|v_{k_0}\|_\infty<C_{k_0}$ by \eqref{norma_limitata} and the validity of 
\eqref{estimate_k0} gives that
$v_{k_0}$ is a positive radial solution of problem \eqref{PROB}
since \eqref{estimate_k0} forces that  $v_{k_0}(r)\le k_0$ for all $r<R$ so that  $g_{k_0}(v_{k_0})=g(v_{k_0})$.

To prove the claim, we suppose by contradiction that $k<\|v_k\|_\infty(\le C_k)$ for all $k\in\mathbb N$ with $v_{k}$ positive solution of truncated problem \eqref{xi0}. 

Using the same change of variables given in \eqref{3.2} and \eqref{defwn}, where $t_k$ and $z_k$ are defined in \eqref{zktk}. So that
from our hypothesis of absurd we have 
\begin{equation} \label{hyp_abs}
t_k\to\infty \quad\text{for} \quad k\to\infty.
\end{equation}
Following the same calculations in Theorem \ref{thm3.1}, we see that for all $k\in\mathbb{N}$ the function $w_k$ is a positive solution of the following problem
\begin{equation} \label{proble_w_k}
  {\footnotesize\begin{cases}
 -\left(\dfrac{y^{N+\alpha-1}|w_k'(y)|^{m-2}w_k'(y)}{\big(a(t_kz_k^{-1}y)+g_k(t_kw_k(y))\big)^\gamma}\right)'=\dfrac{y^{N+\beta-1}}{\big(a(0)+g(k)\big)^\gamma}w_k^p(y),\quad y<\dfrac{z_kR}{t_k},
  \\w_k'(0)=0, \qquad w_k(0)=1, \qquad w_k\biggl(\dfrac{Rz_k}{t_k}\biggr)=0,
   \end{cases}}
\end{equation}
where we have replaced $n$ with $k$ and $\xi=0$ in problem \eqref{prob_w_n}. We can see that 
$w_k'<0$ since $w_k$ is a positive solution of \eqref{proble_w_k} and 
\begin{equation} \label{normaunit}
\|w_k\|_\infty=w_k(0)=\frac{v_k(0)}{\displaystyle\\|v_k\|_\infty}=1    
\end{equation}
for all $k\in\mathbb N$.

Moreover we note that \eqref{R+} still holds 
and, following Theorem \ref{thm3.1}, also  \eqref{limitatezza_derivata} is valid 
for all $k\in\mathbb N$ and $y\in[0,z_kR/t_k)$, from which and since $\beta-\alpha+1>0$ by  $(H_0)$, we get that $w_k'$ is uniformly bounded in compact intervals. For any $\bar R$ positive number then $\bar R<Rz_k/t_k$ for $k$ large and considering the restriction of $w_k$ to $[0,\bar R]$, which we will still denote $w_k$, we can consider that there exists a constant $C(\bar R)>0$ so that 
$$|w_k'(y)|\le C(\bar R), \quad \text{for all} \;\; k\in\mathbb{N}\;\; \text{and}\;\; y\in[0,\bar R],$$
then the sequence $(w_k)_k$ is equilipschitz or equivalently equicontinuous. Furthermore, as we have observed in \eqref{normaunit}, then  $\|w_k\|_\infty=1$ for all $k\in\mathbb N$ then the sequence $(w_k)_k$ is also uniformly bounded in $[0,\bar R]$. By Ascoli Arz\`ela's Theorem, $(w_k)_k$ contains a subsequence converging uniformly  (which we still denote by $(w_k)_k$), namely 
\begin{equation}  \label{conve_w_k}
 w_k\to w \quad\text{in}\;\;C[0,\bar R],  
\end{equation}
furthermore, $\bar R$  can be chosen arbitrary in $\mathbb R^+$ because of $z_kR/t_k\to\infty$ for $k\to\infty$. Consequently, $w$ is well defined in all $\mathbb R^+$ and we immediately get
$\displaystyle\lim_{k\to\infty}w_k(y)=w(y)$ for all $y\in\mathbb R^+$, so that $\displaystyle\lim_{y\to\infty}w(y)=0$ by 
$w_k\bigl(Rz_k/t_k\bigr)=0$ and the validity of \eqref{R+}.
In particular, arguing as we have done to get \eqref{wpositiva}, we obtain that 
\begin{equation}\label{wpositivaR+}
    w(y)>0 \quad\text{for all }\,\,y\in[0,\infty)
\end{equation}
We observe, that since $t_k>k$ for all $k\in\mathbb{N}$ by contradiction, being $t_k=\|v_k\|_\infty$, we have $0<k/t_k<1$. Then, there exist $\ell \in[0,1]$ and a subsequence, which we denote again $(k/t_k)_k$, such that $$\frac{k}{t_k}\to \ell \quad \text{as} \quad k\to\infty.$$
We note that, since $w_k$ is a positive and decreasing solution of problem \eqref{proble_w_k} for all $k\in\mathbb N$, we have that $0< w_k(y)<1$ for all $y\in(0,Rz_k/t_k)$ and $k\in\mathbb N$. In particular, for all $k\in\mathbb N$ there exists $s_k\in(0,Rz_k/t_k)$ such that 
$$w_k(s_k)=\frac{k}{t_k}.$$
We generate a sequence $(s_k)_k\in \mathbb R^+$ which we can assume, without loss of generality, monotone.
We observe that if $s<s_k$, since $w_k$ is a decreasing function for all $k\in\mathbb{N}$, then we have $w_k(s)>w_k(s_k)=k/t_k$ that gives $t_kw_k(s)>k$ for all $s<s_k$. Consequently, by the definition of $g_k$, we get
\begin{equation}\label{=g(k)}
    g_k(t_kw_k(s))=g(k) \quad\text{for all}\,\,s<s_k.
\end{equation}
Now we analyze the limit problem associated with problem \eqref{proble_w_k} by dividing the discussion in three cases: $\ell=0$, $\ell=1$ and $0<\ell<1$.
\begin{enumerate}
\item [(I)] \quad If $\ell=0$, that is $k/t_k\to0$ as $k\to\infty$, we can see that $s_k\to\infty$ as $k\to\infty$. 
Indeed, as we have observed $w_k(y)\to w(y)$ in $\mathbb R^+$. Then, if we suppose by contradiction that $s_k\to\bar s$ for $k\to\infty$, with $\bar s\in\mathbb R^+$ we get 
\begin{equation}\label{limite=0}
0=\lim_{k\to\infty}w_k(s_k)=w(\bar s). 
\end{equation}
On the other hand, \eqref{wpositivaR+} contradicts \eqref{limite=0} then necessarily $\bar s=\infty$.
%on the other hand since $w_k(z_kR/t_k)=0$ then
%$0=\displaystyle\lim_{k\to\infty}w_k(z_kR/t_k)=w(\infty)$.
Thus, for any $\bar R>0$ there is $k_{\bar R}\in\mathbb N$ such that 
  \begin{equation} \label{s_k_bound}
 \bar R< s_k<\frac {z_k R}{t_k}
\quad\text{for all}\;\;k\ge k_{\bar R}.
  \end{equation} 
Now, integrating the equation in \eqref{proble_w_k} from $0$ to $s\in[0,\bar R]$, we get that 
$$|w_k'(s)|^{m-1}=\biggl(\frac{a(t_kz_k^{-1}s)+g_k(t_kw_k(s))}{a(0)+g(k)}\biggr)^\gamma\frac{\int_{0}^{s}\tau^{N+\beta-1}w_k^p(\tau)d\tau}{s^{N+\alpha-1}},$$
then, elevating both members by $1/(m-1)$  we have
\begin{equation} \label{integrazione1}
-w_k'(s)=\biggl(\frac{a(t_kz_k^{-1}s)+g_k(t_kw_k(s))}{a(0)+g(k)}\biggr)^{\frac{\gamma}{m-1}}\biggl[\frac{\int_{0}^{s}\tau^{N+\beta-1}w_k^p(\tau)d\tau}{s^{N+\alpha-1}} \biggr]^{\frac{1}{m-1}}.
\end{equation}
Furthermore, we have that for all $s\in[0,\bar R]$, that implies $s<s_k$ for all $k\ge k_{\bar R}$  by \eqref{s_k_bound}, \eqref{=g(k)} and using also \eqref{limite_frac}, equality \eqref{integrazione1} becomes
\begin{equation} \label{integ_1_kgrande}
-w_k'(s)=\bigl(1+o(1)\bigr)^{\frac{1}{m-1}}\biggl[\frac{1}{s^{N+\alpha-1}} \int_{0}^{s}\tau^{N+\beta-1}w_k^p(\tau)d\tau\biggr]^{\frac{1}{m-1}},
\end{equation}
for all $s\in[0,\bar R]$ and $k\ge k_{\bar R}$.
A further  integration of \eqref{integ_1_kgrande} from $0$ to $y\in[0,\bar R]$ yields for $k\ge k_{\bar R}$
\begin{equation} \label{w_k_equazione}
1-w_k(y)=\bigl(1+o(1)\bigr)^{\frac{1}{m-1}}\int_{0}^{y}\biggl[\frac{1}{s^{N+\alpha-1}}\int_{0}^{s}\tau^{N+\beta-1}w_k^p(\tau)d\tau\biggr]^{\frac{1}{m-1}}ds.
\end{equation}
Now passing to the limit for $k\to\infty$ in \eqref{w_k_equazione}, by using \eqref{conve_w_k} and Lebesgue's dominated Theorem we obtain that $w$ satisfies the following integral equation
\begin{equation} \label{eq_integrale_wk}
1-w(y)=\int_{0}^{y}\biggl[\frac{1}{s^{N+\alpha-1}}\int_{0}^{s}\tau^{N+\beta-1}w^p(\tau)d\tau\biggr]^{\frac{1}{m-1}}ds.
\end{equation}
Now arguing as in \eqref{w'=0} we obtain $w'(0)=0$. In particular, being $(w_k)_k$ a sequence of positive and decreasing functions  then its uniform limit $w$, satisfying \eqref{eq_integrale_wk} and \eqref{wpositivaR+},  is a positive and decreasing solution in $[0,\bar R]$ of the following problem
$$
\begin{cases}
  -\big(y^{N+\alpha-1}|w'(y)|^{m-2}w'(y)\big)'=y^{N+\beta-1}w^p(y), \quad y\in(0,\bar R),
  \\w(0)=1, \qquad w'(0)=0.
\end{cases}
$$  
Being $\bar R$ arbitrary, as already done below formula \eqref{problema_w}, 
 $w$ is a positive solution of \eqref{p2}, contradicting  Theorem \ref{teorema1} applied with $\wp=p$ and $\lambda=1$, so that the case $\ell=0$ cannot occur.
  \item [(II)] \quad If $l=1$, that is $k/t_k\to1$ as $k\to\infty$, we claim that $s_k\to0$ as $k\to\infty$. Indeed, integrating \eqref{integrazione1} from $0$ to $s_k\in(0,z_kR/t_k)$ and using \eqref{=g(k)} we arrive to
 \begin{equation} \label{w_k_equ2.1}
 1-\frac{k}{t_k}=\int_{0}^{s_k}\biggl[\biggl(\frac{a(t_kz_k^{-1}s)+g(k)}{a(0)+g(k)}\biggr)^\gamma\cdot\frac{\int_{0}^{s}\tau^{N+\beta-1}w_k^p(\tau)d\tau}{s^{N+\alpha-1}}\biggr]^{\frac{1}{m-1}}ds.
\end{equation}
By \eqref{limite_frac} and the fact that $w_k$ is a decreasing  function for all $k\in\mathbb N$, so that $w_k(\tau)\ge w_k(s_k)$ for $\tau\le s_k$, from \eqref{w_k_equ2.1} we get 
\begin{equation} \label{w_k_equ2.2}
\begin{aligned}
 1-\frac{k}{t_k}&\ge(1+o(1))^{\frac{1}{m-1}}\biggl(\frac{k}{t_k}\biggr)^{\frac{p}{m-1}}\frac{1}{(N+\beta)^{\frac{1}{m-1}}}\int_0^{s_k}s^\frac{\beta-\alpha+1}{m-1}ds\\
 &=C\bigl(1+o(1)\bigr)^{\frac{1}{m-1}}\biggl(\frac{k}{t_k}\biggr)^{\frac{p}{m-1}}s_k^{\frac{\beta-\alpha+m}{m-1}}\ge 0,
 \end{aligned}
\end{equation}
where $C=(m-1)(N+\beta)^{1/(1-m)}/(\beta-\alpha+m)>0$.
 Now passing to the limit in \eqref{w_k_equ2.2} as $k\to\infty$ and using $k/t_k\to1$ as $k\to\infty$ we obtain that $s_k\to0$ as $k\to\infty$.
 
Consequently for all $y\in(0,\bar R]$ if $k$ is sufficiently large, then $s_k<y$, so that  we can integrate \eqref{integrazione1} from $s_k$ to $y\in[s_k,\bar R]$ obtaining
\begin{equation} \label{w_k_equazione2}
\begin{aligned}
\int_{s_k}^{y}\biggl[\biggl(\frac{a(t_kz_k^{-1}s)+g_k(t_kw_k(s))}{a(0)+g(k)}\biggr)^\gamma&\cdot\frac{\int_{0}^{s}\tau^{N+\beta-1}w_k^p(\tau)d\tau}{s^{N+\alpha-1}}\biggr]^{\frac{1}{m-1}}ds=\frac{k}{t_k}-w_k(y). 
\end{aligned}
\end{equation}
Using \eqref{wpositivaR+} and \eqref{hyp_abs}, given by hypothesis of absurd, we have 
\begin{equation} \label{limite1}
t_kw_k(s)\to\infty \quad\text{for}\quad k\to\infty
\end{equation}
for $s\in[0,\bar R]$ fixed. Then, since \eqref{limite1} and $(H_2)$ we obtain 
\begin{equation} \label{limite2}
g(t_kw_k(s))\sim t_kw_k(s) \quad\text{as}\quad k\to\infty
\end{equation}
for $s\in[0,\bar R]$ fixed. By \eqref{a(0)},  \eqref{limite2} and $t_k\sim k$ as $k\to\infty$, since we are in case $(ii)$, we have that
\begin{equation}\label{comp_asint}
\frac{a(t_kz_k^{-1}s)+g_k(t_kw_k(s))}{a(0)+g(k)}\sim\frac{a(0)+t_kw_k(s)}{a(0)+g(k)}\sim\frac{a(0)+kw_k(s)}{a(0)+k}\sim w_k(s),
\end{equation}
for $k\to\infty$ and for $s\in[0,\bar R]$ fixed, using \eqref{wpositivaR+}.

Now passing to the limit for $k\to\infty$ in \eqref{w_k_equazione2}, by using \eqref{conve_w_k}, Lebesgue's dominated Theorem and \eqref{comp_asint} we obtain that $w$ satisfies the following integral equation
\begin{equation} \label{eq_integrale_wk2}
1-w(y)=\int_{0}^{y}\biggl[\frac{w(s)^\gamma}{s^{N+\alpha-1}}\int_{0}^{s}\tau^{N+\beta-1}w^p(\tau)d\tau\biggr]^{\frac{1}{m-1}}ds.
\end{equation}
Now arguing as in \eqref{w'=0} we obtain $w'(0)=0$. In particular, being $(w_k)_k$ a sequence of positive and decreasing functions  then its uniform limit $w$, satisfying \eqref{eq_integrale_wk2},  is a positive and decreasing solution in $[0,\bar R]$ of the following problem
$$
\begin{cases}
  -\biggl(\dfrac{y^{N+\alpha-1}|w'(y)|^{m-2}w'(y)}{w^\gamma(y)}\biggr)'=y^{N+\beta-1}w^p(y), \quad y\in(0,\bar R),
  \\w(0)=1, \qquad w'(0)=0.
\end{cases}
$$
 %Following the same argument using at the end of the proof of Theorem \ref{thm3.1} we can extended $w$ to the entire $\mathbb{R}^+$.
 Let us consider the following change of variables
 $$u(y)=w^{1-\frac{\gamma}{m-1}}(y),$$ 
  then  we have that $u$ is a positive solution in $[0,\bar R]$ of the following problem
$$
\begin{cases}
    -\bigl(y^{N+\alpha-1}|u'(y)|^{m-2}u'(y)\bigr)'\!=\bigl(\!1-\frac{\gamma}{m-1}\!\bigr)^{m-1} \!y^{N+\beta-1}\!u^{\frac{p(m-1)}{m-1-\gamma}}(y), \,\,y\in(0,\bar R)\\
    u(0)=1, \qquad u'(0)=0.
\end{cases}
$$
This is problem \eqref{p2} with 
\begin{equation}\label{lambda&p}
\lambda=\biggl(1-\dfrac{\gamma}{m-1}\biggr)^{m-1} \quad\text{and}\quad \wp=\frac{p(m-1)}{m-1-\gamma},    
\end{equation}
such that $\wp$ satisfies $(H_1)'$ since $p$ satisfies $(H_1)$. 
Following the same argument using at the end of the proof of Theorem \ref{thm3.1}, we get that $u$ can be extended to the entire $\mathbb{R}^+$, obtaining a positive solution of \eqref{p2} with $\lambda$ and $\wp$ defined in \eqref{lambda&p}.
This contradicts Theorem \ref{teorema1}, thus also the case $\ell=1$ cannot occur.

\item [(III)] \quad If $0<\ell<1$, we can see that $(s_k)_k$ is bounded. We first observe that, if we suppose that $s_k\to\infty$ as $k\to\infty$, passing to the limit for $k\to\infty$ in \eqref{w_k_equ2.2}, that holds for every limit of $k/t_k$, we get
$$
1-\ell\ge C\bigl(1+o(1)\bigr)^{\frac{\gamma}{m-1}}\ell^{\frac{p}{m-1}}\lim_{k\to\infty}s_k^{\frac{\beta-\alpha+m}{m-1}}=\infty 
$$
since $\beta-\alpha+m>0$ by $(H_0)$ and $m>1$, that contradicts $\ell\in(0,1)$. Then the sequence $(s_k)_k$ is bounded. Then by compactness, there is $s_0\in\mathbb R^+$ and a subsequence, which we denote again by $(s_k)_k$, such that $s_k\to s_0$.
Using \eqref{=g(k)} we get that $w_k$ satisfies the integral equation \eqref{w_k_equazione} in $(0,s_k)$.
Proceeding as in case $(i)$, we conclude that $w$ is a positive solution in $[0,s_0]$ of the following problem
\begin{equation} \label{problema_w4}
\begin{cases}
  -\big(y^{N+\alpha-1}|w'(y)|^{m-2}w'(y)\big)'=y^{N+\beta-1}w^p(y), \quad y\in(0,s_0),
  \\w(0)=1, \qquad w'(0)=0.
\end{cases}
\end{equation}
On the other hand, for each $s\in(s_k,\infty)$ fixed, we can see that 
\begin{equation}\label{beha}
t_kw_k(s)\sim \frac{k}{\ell}w_k(s) \quad \text{as} \quad k\to\infty
\end{equation}
since $k/t_k\to \ell$. Moreover, by \eqref{wpositivaR+} we obtain that $t_kw_k(s)\to\infty$ for $k\to\infty$ and $s\in(s_k,\infty)$ fixed. Using \eqref{a(0)}, \eqref{beha} and $(H_2)$ we get 
\begin{equation} \label{comp_asint_2}
\frac{a(t_kz_k^{-1}s)+g_k(t_kw_k(s))}{a(0)+g(k)}\sim \frac{a(0)+kw_k(s)/\ell}{a(0)+k}\sim \frac{1}{\ell}w_k(s),
\end{equation}
for $s\in(s_k,\infty)$ fixed and $k\to\infty$.
From \eqref{w_k_equazione2}, passing to the limit for $k\to\infty$ and using \eqref{conve_w_k}, Lebesgue's dominated Theorem and \eqref{comp_asint_2} we obtain that $w$ satisfies the following integral equation
$$
\ell-w(y)=\int_{s_0}^{y}\biggl[\frac{(l^{-1}w_k(s))^\gamma}{s^{N+\alpha-1}}\int_{0}^{s}\tau^{N+\beta-1}w_k^p(\tau)d\tau\biggr]^{\frac{1}{m-1}}ds, 
$$
for all $y\in(s_0,\infty)$. Thus, since $w_k(s_k)\to \ell$ but also $w_k(s_k)\to w(s_0)$, we have that $w$ is a positive solution in $(s_0,\infty)$ of the initial value problem 
\begin{equation} \label{prob_2parte}
\begin{cases}
      -\left(\frac{y^{N+\alpha-1}|w'(y)|^{m-2}w'(y)}{\ell^{-\gamma} w^\gamma(y)}\right)'= y^{N+\beta-1}w^{p}(y),& y\in (s_0,\infty), \\w(s_0)=\ell.
   \end{cases}
\end{equation}
Finally, from \eqref{problema_w4} and \eqref{prob_2parte}, we get that $w$ is a positive solution of \eqref{P} with $d=1/\ell>1$ since $l\in(0,1)$. This contradicts Theorem \ref{teorema2}.
\end{enumerate}
Consequently, condition \eqref{estimate_k0} holds and the proof of the theorem is concluded.
\end{proof}

\section{Nonexistence results} \label{cap3}
In this section we develop the proof of the two nonexistence results  of positive radial solutions of \eqref{PROB}, given by  Theorem \ref{th_1.2} and Theorem \ref{th_1.3} stated in the Introduction.

The proof of the first nonexistence result, that is Theorem \ref{th_1.2}, is obtained via 
of a Poho\v zaev-Pucci-Serrin type identity
for positive radial solutions. Beyond Poho\v zaev's identity in \cite{P} (cfr. \cite{ZN}),  a pioneering radial identity for quasilinear problems  can be found \cite{NS2}, and \cite{NS},  where Ni and Serrin considered positive radial solutions of \eqref{nsprob}. Our new radial identity is rather delicate since it involves a new exponent $\sigma$ to be chosen properly, which makes calculations quite cumbersome. 
\begin{proposition}
Assume $(H_0)$.
%$(H_2)'$ and $(H_3)'$.
Let $v\in C^1[0,R]$  be  positive radial solution of problem \eqref{PROB} with $a$ and $g$ of class $C^1$, then for a constant
\begin{equation}\label{sigma_bounds}-(N+\alpha-m)<\sigma\le m-1\end{equation}
the following identity holds
\begin{equation}\label{final}
\begin{aligned}
\biggl(\sigma-&m+2-\frac{N+\alpha-m+1+\sigma}{m}+\frac{N+\beta+1+\sigma-m}{p+1}\biggr)\!\!\int_{0}^{R}\frac{r^{N+\alpha-m+\sigma}|v'(r)|^m}{\big(a(r)+g(v(r))\big)^\gamma}dr\\
&=\frac{m-1}{m}\frac{R^{N+\alpha-m+1+\sigma}|v'(R)|^m}{\big(a(R)+g(0)\big)^\gamma}-\frac{\gamma}{m}\int_{0}^{R}\frac{r^{N+\alpha-m+1+\sigma}|v'(r)|^ma'(r)}{\big(a(r)+g(v(r))\big)^{\gamma+1}}dr\\
&\quad-\frac{\gamma}{m}\int_{0}^{R}\frac{r^{N+\alpha-m+1+\sigma}|v'(r)|^{m}v'(r)g'(v(r))}{\big(a(r)+g(v(r))\big)^{\gamma+1}}dr\\
&\quad{{\,-\frac{N+\beta+1+\sigma-m}{p+1}(\sigma-m+1)\int_0^R\dfrac{r^{N+\alpha+\sigma-m-1}|v'(r)|^{m-2}v'(r)v(r)}{\bigl(a(r)+g(v(r)))^\gamma}dr}}.
\end{aligned}\end{equation}
\end{proposition}
\begin{proof}
Let $v\in C^1[0,R]$  be a  positive radial solution of problem \eqref{PROB},
or equivalently a positive solution of \eqref{probrad}. By  \eqref{derivatasegno} and \eqref{reg2}, $v$ is a strictly decreasing function in $[0,R]$  with $v\in C^2(0,R)$.

Multiplying the equation in \eqref{probrad} by $r^{\sigma-m+2} v'(r)$, with
$\sigma$ as in \eqref{sigma_bounds}, and integrating from $0$ to $R$, we find
\begin{equation} \label{6.1}
    -\int_{0}^{R}\biggl(\frac{r^{N+\alpha-1}|v'(r)|^{m-2}v'(r)}{\big(a(r)+g(v(r))\big)^\gamma}\biggr)'r^{\sigma-m+2} v'(r)dr=\int_{0}^{R}r^{N+\beta+\sigma-m+1}v^p(r)v'(r)dr.
\end{equation}
Now we rewrite the expression \eqref{6.1} as
\begin{equation}\label{A=B}
    \mathcal{L}=\mathcal{R}.
\end{equation}
The term on the left hand side, integrating by parts and using that $N+\alpha+\sigma-m+1>0$ by \eqref{sigma_bounds}  and  $a(R)\ge c_1>0$ by $(H_3)$,  becomes
$$\begin{aligned}
\mathcal{L}=&-\frac{R^{N+\alpha+\sigma-m+1}|v'(R)|^m}{\big(a(R)+g(0)\big)^\gamma}+(\sigma-m+2)\int_{0}^{R}\frac{r^{N+\alpha+\sigma-m}|v'(r)|^m}{\big(a(r)+g(v(r))\big)^\gamma}dr\\&
\int_{0}^{R}\frac{r^{N+\alpha+\sigma-m+1}|v'(r)|^{m-2}v'(r)v''(r)}{\big(a(r)+g(v(r))\big)^\gamma}dr.
\end{aligned}$$
In particular,
\begin{equation}\label{A}
\begin{aligned}
&\mathcal{L}=-\frac{m-1}{m}\frac{R^{N+\alpha-m+1+\sigma}|v'(R)|^m}{\big(a(R)+g(0)\big)^\gamma}\\
&\qquad+\biggl(\sigma-m+2-\frac{N+\alpha-m+1+\sigma}{m}\biggr)\int_{0}^{R}\frac{r^{N+\alpha-m+\sigma}|v'(r)|^m}{\big(a(r)+g(v(r))\big)^\gamma}dr\\
&\qquad+\frac{\gamma}{m}\int_{0}^{R}\frac{r^{N+\alpha-m+1+\sigma}|v'(r)|^ma'(r)}{\big(a(r)+g(v(r))\big)^{\gamma+1}}dr\\
&\qquad+\frac{\gamma}{m}\int_{0}^{R}\frac{r^{N+\alpha-m+1+\sigma}|v'(r)|^{m}v'(r)g'(v(r))}{\big(a(r)+g(v(r))\big)^{\gamma+1}}dr.
\end{aligned}\end{equation}
On the other hand, concerning , integrating by parts $\mathcal{R}$ and using $v(R)=0$, we get
\begin{equation}\label{6.3}
\begin{aligned}
\mathcal{R}&=\int_{0}^{R}r^{N+\beta+\sigma-m+1}v^p(r)v'(r)dr=\frac{1}{p+1}\int_{0}^{R}r^{N+\beta+\sigma-m+1}\big(v^{p+1}(r)\big)'dr\\
&=-\frac{N+\beta+\sigma-m+1}{p+1}\int_{0}^{R}r^{N+\beta+\sigma-m}v^{p+1}(r)dr,
\end{aligned}\end{equation}
since, by \eqref{sigma_bounds} and $\beta-\alpha+1>0$,
it trivially holds
\begin{equation}\label{sigma_beta}N+\beta+\sigma-m+1> N+\alpha+\sigma-m>0.\end{equation}
Similarly, multiplying the equation in \eqref{probrad} by $r^{\sigma-m+1}v(r)$ and  integrating from $0$ to $R$ we get 
$$
\int_0^R\dfrac{r^{N+\alpha-1}|v'(r)|^{m-2}v'(r)}{\bigl(a(r)+g(v(r)))^\gamma}\bigl[(\sigma-m+1)r^{\sigma-m}v(r)+r^{\sigma-m+1}v'(r)\bigr]dr
=\int_0^R r^{N+\beta-m+\sigma}v^{p+1}(r)$$
since $v'(0)=v(R)=0$
and  $N+\alpha+\sigma-m>0$ by \eqref{sigma_bounds}.
Consequently
\begin{equation}\label{6.4}
\begin{aligned}(\sigma-m+1)\int_0^R&\dfrac{r^{N+\alpha+\sigma-m-1}|v'(r)|^{m-2}v'(r)v(r)}{\bigl(a(r)+g(v(r)))^\gamma}dr+
\int_0^R\dfrac{r^{N+\alpha+\sigma-m}|v'(r)|^{m}}{\bigl(a(r)+g(v(r)))^\gamma}
dr\\
&\hspace{1cm}=\int_0^R r^{N+\beta-m+\sigma}v^{p+1}(r)dr.   
\end{aligned}\end{equation}
Combining \eqref{6.3} and \eqref{6.4}, we obtain
\begin{equation} \label{B}\begin{aligned}
\mathcal{R}=&-\frac{N+\beta+1+\sigma-m}{p+1}(\sigma-m+1)\int_0^R\dfrac{r^{N+\alpha+\sigma-m-1}|v'(r)|^{m-2}v'(r)v(r)}{\bigl(a(r)+g(v(r)))^\gamma}dr\\
&-\frac{N+\beta+1+\sigma-m}{p+1}
\int_0^R\dfrac{r^{N+\alpha+\sigma-m}|v'(r)|^{m}}{\bigl(a(r)+g(v(r)))^\gamma}
dr.
\end{aligned}\end{equation}
Then, using \eqref{A} and \eqref{B} in \eqref{A=B}, the identity\eqref{final} follows at once.
\end{proof}

Now we are ready to prove the main nonexistence result, that is Theorem \ref{th_1.2}, whose statement is given in the Introduction.

\begin{proof}[Proof of Theorem \ref{th_1.2}]
Assume by contradiction that there exists $v\in C^1[0,R]$   positive radial solution of problem \eqref{PROB},
or equivalently positive solution of \eqref{probrad}. By  \eqref{derivatasegno} and \eqref{reg2}, $v$ is a strictly decreasing function in $[0,R]$  with $v\in C^2(0,R)$.

Now, since $a'\le0$ and $g'\ge0$, by $(H_2)'$ and $(H_3)'$, and using monotonicity of $v$, the first three terms on the right hand side of the identity \eqref{final} are positive, yielding
\begin{equation}\label{contrad}c_1\int_{0}^{R}\frac{r^{N+\alpha-m+\sigma}|v'(r)|^m}{\big(a(r)+g(v(r))\big)^\gamma}dr-c_2\int_0^R\dfrac{r^{N+\alpha+\sigma-m-1}|v'(r)|^{m-1}v(r)}{\bigl(a(r)+g(v(r)))^\gamma}dr>0,\end{equation}
where
$$c_1=\sigma-m+2-\frac{N+\alpha-m+1+\sigma}{m}+\frac{N+\beta+1+\sigma-m}{p+1}$$ and
$$c_2=\frac{N+\beta+1+\sigma-m}{p+1}(m-1-\sigma)$$
and with $\sigma$ satisfying  \eqref{sigma_bounds}.
Furthermore, thanks to \eqref{sigma_beta} and \eqref{sigma_bounds}, it follows  $c_2\ge0$.

We claim that $c_1<0$ if $p+1> m^*_{\alpha,\beta,\gamma}$. In particular, this latter condition is equivalent to 
$$p+1> \dfrac{m(N+\beta-m+1+\sigma)}{N+\alpha-m+(m-1)(m-1-\sigma)},$$
whenever  $\sigma$ is chosen such that 
\begin{equation}\label{cond_for_sigma}
\dfrac{m(N+\beta-m+1+\sigma)}{N+\alpha-m+(m-1)(m-1-\sigma)}
=m^*_{\alpha,\beta,\gamma},
\end{equation}
with $m^*_{\alpha,\beta,\gamma}$ given in $(H_1)$.
Condition \eqref {cond_for_sigma}
yields the explicit expression of $\sigma$, given by
$$
\sigma=\dfrac{m(m-1)^2\bigl[\alpha-\beta +m(N+\beta-1)\bigr]-\gamma[N+\alpha-m+(m-1)^2][m(N+\beta+1)-N-\alpha]}{(m-1)\bigl\{m[\alpha-\beta+m(N+\beta-1)]-\gamma[m(N+\beta+1)-N-\alpha]\bigr\}},
$$
where $m(N+\beta+1)-N-\alpha>0$, by \eqref{first_sign}, while  $(H_0)$ gives
$$\alpha-\beta+m(N+\beta-1)>(m-1)(\beta-\alpha+m)>0.$$
In addition
\begin{equation}\label{gamma_2value}m[\alpha-\beta+m(N+\beta-1)]-\gamma[m(N+\beta+1)-N-\alpha]>0\end{equation}
by the choice of $\gamma$.
Indeed, \eqref{gamma_2value} is equivalent to
$$\gamma<\dfrac{m[\alpha-\beta+m(N+\beta-1)]}{m(N+\beta+1)-N-\alpha},$$
which follows from the upper bound $\Upsilon$ for $\gamma$ given in $(H_1)$, being
$$\dfrac{m[\alpha-\beta+m(N+\beta-1)]}{m(N+\beta+1)-N-\alpha}>\Upsilon,$$
by $(H_0)$.

We are now ready to verify that  $\sigma$ above satisfies \eqref{sigma_bounds}.
First we note that inequality
$\sigma\le m-1$ is equivalent to 
$$\begin{aligned}
\gamma(N+&\alpha-m)[N+\alpha-m(N+\beta+1)]
+\gamma (m-1)^2[N+\alpha-m(N+\beta+1)]\\
&\qquad+m(m-1)^2\bigl[\alpha-\beta +m(N+\beta-1)\bigr]\\
&\le m(m-1)^2[\alpha-\beta+m(N+\beta-1)]+
\gamma (m-1)^2[N+\alpha-m(N+\beta+1)]
\end{aligned}$$
which holds since
$$\gamma(N+\alpha-m)[N+\alpha-m(N+\beta+1)]
<0$$
by $(H_0)$ and \eqref{first_sign}.

To verify that $\sigma>m-N-\alpha$, we need to show that
$$\begin{aligned}\gamma&[N+\alpha-m+(m-1)^2][N+\alpha-m(N+\beta+1)]+m(m-1)^2\bigl[\alpha-\beta +m(N+\beta-1)\bigr]\\
&+(N+\alpha-m)(m-1)\biggl\{m[\alpha-\beta+m(N+\beta-1)]-\gamma[m(N+\beta+1)-N-\alpha]\biggr\}>0
\end{aligned}$$
namely
$$\begin{aligned}
\gamma&[m(N+\beta+1)-N-\alpha]\biggl[N+\alpha-m+(m-1)^2+(N+\alpha-m)(m-1)\biggr]\\
&<m(m-1)^2\bigl[\alpha-\beta +m(N+\beta-1)\bigr]+(N+\alpha-m)(m-1)m[\alpha-\beta+m(N+\beta-1)].
\end{aligned}$$
The previous condition holds if
$$\gamma[m(N+\beta+1)-N-\alpha][m(N+\alpha-1)-m+1]<m(m-1)\bigl[\alpha-\beta +m(N+\beta-1)\bigr](N+\alpha-1)$$
which, by virtue of \eqref{first_sign} and  using that
$$m(N+\alpha-1)-m+1= m(N+\alpha-m)+(m-1)^2>0,$$
is equivalent to  $\gamma<\Upsilon_1$, with $\Upsilon_1$ is defined in \eqref{ups1}.
We need now to compare $\Upsilon$ given in $(H_1)$ and $\Upsilon_1$. Precisely, $\Upsilon>\Upsilon_1$  holds if and only if 
$$\dfrac{\bigl[\alpha-\beta +m(N+\beta-1)\bigr](N+\alpha-1)}{m(N+\alpha-1)-m+1}<\beta-\alpha+m$$
that is 
$$(N+\alpha-1)\bigl[\alpha-\beta+m(N+\beta-1)-m(\beta-\alpha+m)\bigr]+(m-1)(\beta-\alpha+m)<0$$
which gives
%$$(N+\alpha-1)[\alpha-\beta-m+m(N+\alpha-m)]+(m-1)(\bet%a-\alpha+m)<0$$
%or equivalently
$$(N+\alpha-m)\bigl[m(N+\alpha-1)+\alpha-\beta-m\bigr]<0$$
yielding
\begin{equation}\label{sigma_ge}
m(N+\alpha-1)-m+\alpha-\beta<0,
\end{equation}
which is exactly the  inequality in $(ii)$ in the statement of Theorem \ref{th_1.2}. 
Consequently, only if \eqref{sigma_ge} holds we need to restrict the range of $\gamma$ in $(H_1)$, by requiring $\gamma<\Upsilon_1$ to obtain a suitable $\sigma$ satisfying \eqref{sigma_bounds}.
In turn, the proof of the theorem is so concluded, since the required contradiction is obtained in \eqref{contrad} being
$c_1<0$ and $c_2\ge0$.
\end{proof}

\begin{remark}
We point out that, from $N+\alpha>m$, condition
$m(N+\alpha-1)-m+\alpha-\beta\ge0$
is implied by
\begin{equation}\label{sigma_pos}
(m-2)m\ge\beta-\alpha,
\end{equation}
which is the condition yielding $\sigma>0$. This latter, by previous calculations,  is equivalent to  
$$\gamma<\dfrac{m(m-1)^2\bigl[\alpha-\beta +m(N+\beta-1)\bigr]}{(N+\alpha-m+(m-1)^2)[m(N+\beta+1)-N-\alpha]}:=\Upsilon_2$$
which follows from the upper bound $\Upsilon$ for $\gamma$ when \eqref{sigma_pos} holds, being $\Upsilon_2>\Upsilon$.
\end{remark}

Using properties of positive solutions of \eqref{probrad}, we obtain a second nonexistence result, Theorem \ref{th_1.3} whose statement is given in the Introduction, which investigates the complementary condition on the sign of $\beta-\alpha+m$, respect to that assumed in Theorem~\ref{th_1.2}.

\begin{proof}[Proof of Theorem \ref{th_1.3}]
Assume by contradiction that there exists $v\in C^1(0,R)\cap C[0,R]$ a positive solution of problem \eqref{probrad}. We have that $v$ verifies the following integral equation
$$v(r)=\int_{r}^{R}\biggl[\frac{\big(a(s)+g(v(s))\big)^\gamma}{s^{N+\alpha-1}}\int_{0}^{s}\tau^{N+\beta-1}v^p(\tau)d\tau\biggr]^{\frac{1}{m-1}}ds.$$
Since $v$ is a positive and decreasing function by \eqref{derivatasegno} with $v(R)=0$ we have that, given $\varepsilon\in(0,1)$, there exists $r_0\in(0,R/2)$ such that
\begin{equation} \label{r_0}
  v(r_0)>\varepsilon.  
\end{equation}
Then, by \eqref{r_0}, since $g$ is nonnegative and $v$ decreasing, we have, for all $r<r_0$,
\begin{equation}\label{v(r)}
\begin{aligned}
v(r)&\ge c_1^\gamma\int_{r}^{r_0}\biggl[\frac{1}{s^{N+\alpha-1}}\int_{0}^{s}\tau^{N+\beta-1}v^p(\tau)d\tau\biggr]^{\frac{1}{m-1}}ds\\
&\ge c_1^\gamma\biggl(\frac{v^p(r_0)}{N+\beta}\biggr)^{\frac{1}{m-1}}\int_{r}^{r_0}s^{\frac{\beta-\alpha+1}{m-1}}ds \ge K\int_{r}^{r_0}s^{\frac{\beta-\alpha+1}{m-1}}ds,\end{aligned}
\end{equation}
where
$K=\varepsilon^{p/(m-1)}c_1^\gamma(N+\beta)^{-\frac{1}{m-1}}.$
Since $\beta-\alpha+m\le0$ by assumption, then $\beta-\alpha+1<0$ being $m>1$. Now, we divide the proof in two cases. 

If $\beta-\alpha+m<0$, then \eqref{v(r)} becomes
$$v(r)\ge K\frac{\alpha-\beta-m}{m-1}\big(r^{\frac{\beta-\alpha+m}{m-1}}-r_0^{\frac{\beta-\alpha+m}{m-1}}\big), \quad 0<r<r_0.$$

While, if $\beta-\alpha+m=0$, from inequality \eqref{v(r)} we obtain
$$
v(r)\ge K\int_{r}^{r_0}s^{-1}ds=K\text{ln}\biggl(\frac{r_0}{r}\biggr), \quad 0<r<r_0.
$$
In both cases we get $v(0)=\lim_{r\to0^+}v(r)=\infty$.
This is a contradiction since $v\in C[0,R]$. The proof of the Theorem \ref{th_1.3} is so concluded.
\end{proof}

\section*{Acknowledgments}
R. Filippucci and L. Baldelli are members of the {\em Gruppo Nazionale per
l'Analisi Ma\-te\-ma\-ti\-ca, la Probabilit\`a e le loro Applicazioni}
(GNAMPA) of the {\em Istituto Nazionale di Alta Matematica} (INdAM).
R. Filippucci was partly supported by  {\em Fondo Ricerca di
Base di Ateneo Esercizio} 2017-19 of the University of Perugia, named {\em Problemi con non linearit\`a dipendenti dal gradiente}.
L. Baldelli was partially supported by National Science Centre, Poland (Grant No. 2020/37/B/ST1/02742).

\end{document}